\documentclass[a4paper,11pt]{amsart}

\usepackage[centertags]{amsmath}
\usepackage{amsfonts}
\usepackage{amssymb}
\usepackage{amsthm}
\usepackage{newlfont}
\usepackage{euscript}

\newtheorem{thm}{Theorem}[section]
\newtheorem{lem}{Lemma}[section]
\newtheorem{cor}{Corollary}[section]
\newtheorem{prop}{Proposition}[section]
\theoremstyle{remark}
\newtheorem{rem*}{Remark}[section]
\theoremstyle{definition}
\newtheorem{dfn}{Definition}[section]

 \numberwithin{equation}{section}

\begin{document}
\newcommand{\No}{N}
\renewcommand{\phi}{\varphi}
\newcommand{\1}{1\!\!{\mathrm I}}
\renewcommand{\Re}{{\Bbb R}}
\newcommand{\eps}{\varepsilon}
\newcommand{\kap}{\varkappa}
\newcommand{\vO}{\varOmega}
\newcommand{\bu}{\mathbf{u}}
\newcommand{\bp}{\mathbf{p}}
\newcommand{\bx}{{x_*}}
\newcommand{\bc}{\mathbf{c}}
\newcommand{\bt}{t_*}
\newcommand{\bq}{\mathbf{q}}
\newcommand{\bR}{\mathbf{R}}
\newcommand{\bY}{\mathbf{Y}}
\newcommand{\ax}{\Re^+}
\newcommand{\prt}{\partial}
\newcommand{\Es}{\mathsf{E}}
\newcommand{\Cs}{\mathsf{C}}
\newcommand{\Ps}{\mathsf{P}}
\newcommand{\Xs}{\mathsf{X}}
\newcommand{\Ys}{\mathsf{Y}}
\newcommand{\Ls}{\mathsf{L}}
\newcommand{\Ks}{\mathsf{K}}
\newcommand{\Ff}{{\EuScript F}}
\newcommand{\Xf}{\mathfrak{X}}
\newcommand{\Yf}{{\EuScript Y}}
\newcommand{\Ef}{{\EuScript E}}
\newcommand{\Bf}{\mathfrak{B}}
\newcommand{\Tf}{{\EuScript T}}
\newcommand{\Cf}{{\EuScript C}}
\newcommand{\Zf}{{\EuScript Z}}
\newcommand{\Wf}{{\EuScript W}}
\newcommand{\Df}{{\EuScript D}}
\newcommand{\Nf}{{\EuScript N}}
\newcommand{\Gf}{{\EuScript G}}
\newcommand{\Qf}{{\EuScript Q}}
\newcommand{\Af}{\mathfrak{A}}
\newcommand{\Lf}{{\EuScript L}}
\newcommand{\Jf}{{\EuScript J}}
\newcommand{\Uf}{{\EuScript U}}
\newcommand{\Sf}{{\EuScript S}}
\newcommand{\Pf}{{\EuScript P}}
\newcommand{\Kf}{{\EuScript K}}
\newcommand{\pf}{\Pi_{fin}}
\newcommand{\Kb}{\mathbf{K}}
\newcommand{\ZZ}{{\Bbb Z}}
\newcommand{\NN}{{\Bbb N}}
\newcommand{\QQ}{{\Bbb Q}}
\newcommand{\TT}{{\Bbb T}}
\newcommand{\GG}{{\Bbb G}}
\newcommand{\DD}{{\Bbb D}}
\newcommand{\NNN}{\NN\times \NN}
\newcommand{\demo}{\emph{Proof.} }
\newcommand{\Cd}{C_\bullet\,}
\newcommand{\Span}{\mathrm{span}\,}
\newcommand{\supp}{\mathrm{supp}\,}
\newcommand{\pa}{\prt_\alpha}
\newcommand{\pba}{\prt_{\sba}}
\newcommand{\be}{\begin{equation}}
\newcommand{\ee}{\end{equation}}
\newcommand{\nvar}[1]{\left\|#1\right\|_{var}}
\newcommand{\bmu}{{\mu_{inv}}}
\newcommand{\tov}{\mathop{\longrightarrow}\limits^{var}}
\newcommand{\eqd}{\mathop{=}\limits^{df}}
\newcommand{\indicator}[1]{1\!\!\mathrm I\{#1\}}
\newcommand{\brakes}[1]{\left(#1\right)}
\newcommand{\eqdistr}{\stackrel{\mathrm{d}}{=}}
\newcommand{\eqdef}{\stackrel{\mathrm{def}}{=}}
\newcommand{\limw}{\stackrel{\mathrm{w}}{\rightarrow}}
\newcommand{\lims}{\underset{n\rightarrow\infty}{{\lim\sup}}\,}
\newcommand{\norm}[1]{\left\Vert#1\right\Vert}
\newcommand{\abs}[1]{\left\vert#1\right\vert}
\newcommand{\figbrakes}[1]{\left\{#1\right\}}
\newcommand{\sqrbrakes}[1]{\left[#1\right]}
\newcommand{\mes}[1]{mes\{#1\}}
\newcommand{\Bc}{\mathcal{B}}
\newcommand{\Gc}{\mathcal{G}}
\newcommand{\Fc}{\mathcal{F}}
\newcommand{\Nc}{\mathcal{N}}
\newcommand{\Lc}{\mathcal{L}}

%\selectlanguage{english}

\title[Malliavin calculus for difference approximations of multidimensional diffusions]
    {Malliavin calculus for difference approximations of multidimensional diffusions:
    truncated local limit theorem}

%    Information for first author
\author{Alexey M. Kulik}
\address{Institute of Mathematics,
Ukrai\-ni\-an National Academy of Sciences, 3, Tereshchenkivska
Str., Kyiv 01601, Ukraine}
 \email{kulik@imath.kiev.ua}

\begin{abstract}
For a difference  approximations of multidimensional diffusion,
the truncated local limit theorem is proved. Under very mild
conditions on the distribution of the difference terms, this
theorem provides that the transition probabilities of these
approximations, after truncation of some asymptotically negligible
terms, possess a densities that converge uniformly to the
transition probability density for the limiting diffusion and
satisfy a uniform diffusion-type estimates. The proof is based on
the new version of the Malliavin calculus for the product of
finite family of measures, that may contain non-trivial singular
components. An applications for uniform estimates for mixing and
convergence rates for difference approximations to SDE's and for
convergence of difference approximations for local times of
multidimensional diffusions are given.

\end{abstract} 


\thanks{The research was partially supported by the Ministry
of Education and Science of Ukraine, project N GP/F26/0106}

 \maketitle

\section*{Introduction}

Consider diffusion process $X$ in $\Re^d$ defined by an SDE
\be\label{00}X(t)=X(0)+\int_0^ta(X(s))\,ds+\int_0^tb(X(s))\,
dW(s), \quad t\in \ax, \ee and a sequence of processes $X_n, n\geq
1$, with their values  at the time moments ${k\over n}, k\in \NN$
defined  by a difference relation \be\label{01}X_n\brakes{k\over
n}=X_n\brakes{k-1\over n}+a\brakes{X_n\brakes{k-1\over n}}\cdot
{1\over n} +b\brakes{X_n\brakes{k-1\over n}}\cdot {\xi_k\over\sqrt
n}, \ee and, at all the other time moments, defined in a
piece-wise linear way:
 \be\label{02} X_n(t)=X_n\brakes{k-1\over
n}+(nt-k+1)\left[X_n\brakes{k\over n}-X_n\brakes{k-1\over n}
\right],\quad t\in\left[{k-1\over n},{k\over n}\right). \ee Here
and below,  $W$ is a Wiener process valued in $\Re^d$, $\{\xi_k\}$
is a sequence of i.i.d. random vectors in $\Re^d$, that belong to
the domain of attraction of the normal law, are centered and have
the identity for covariance matrix. Under standard assumptions on
the coefficients of the equations (\ref{00}), (\ref{01}) (local
Lipschitz condition and linear growth condition), the
distributions of the processes $X_n$ in  $C(\ax, \Re^d)$ with the
given initial value $X_n(0)=x$ converge weakly to the distribution
of the process $X$ with $X(0)=x$ (\cite{skor_as_m}). Thus, it is
natural to call the sequence $\{X_n\}$ \emph{the difference
approximation} for the diffusion $X$.

Consider the transition probabilities for the processes $X$,$X_n$:
$$P_{x,t}(dy)
\equiv P(X(t)\in dy|X(0)=x),\quad P_{x,t}^n(dy)\equiv P(X_n(t)\in
dy|X_n(0)=x),\quad t>0, x\in\Re^d.$$  It is well known
(\cite{IKO}) that if the coefficients $a,b$ are H\"older
continuous and bounded and the matrix $b\cdot b^*$ is uniformly
non-degenerate, then $P_{x,t}(dy)=p_t(x,y)\,dy$. The function
$\{p_{t}(x,y), t\in\ax, x,y\in\Re^d\}$ (the \emph{transition
probability density} for $X$) possesses the  estimate

\be\label{A2} p_t(x,y)\leq C(T)\,t^{-{d\over
2}}\exp\brakes{-{\gamma\|y-x\|^2\over t}}, \quad t\leq T,
x,y\in\Re^d.
\ee

The general question, that motivates the present paper, is whether
any (more or less restrictive) conditions can be imposed on the
coefficients $a,b$ and the distribution of $\xi_k$ in order to
provide that $P^n_{x,t}(dy)=p^n_t(x,y)\,dy$ for $n$ large enough,
the densities $p^n$ possess an estimate analogous to (\ref{A2})
and $p^n$ converge to $p$ in an appropriate way. Such a question
both is interesting by itself and has its origin in the numerous
applications, such as nonparametric estimation problems in time
series analysis and diffusion models (see the discussion in the
Introduction to \cite{kon_mammen}), the uniform bounds for the
mixing coefficients of the difference approximations to SDE's (see
\cite{Klo_Ver} and subsection 4.1 below), the difference
approximation for local times of multidimensional diffusions (see
\cite{kulik_da} and subsection 4.2 below).

In the current paper, we consider the  question exposed above in a
slightly modified setting. For the distributions $P^n$, we prove
the result that we call \emph{the truncated local limit theorem}.
Let us explain this term. We show that the kernel $P^n$ can be
decomposed into the sum $P^n=Q^n+R^n$ in such a way that both
$Q^n$ and $R^n$ are a non-negative kernels and

(i) for $Q^n$, its density $q^n$ exists, satisfies an analogue of
 (\ref{A2}) and converges to $p$;

(ii) for $R^n$, its total mass can be estimated explicitly and
converges to $0$.

The kernel  $Q^n$ represents \emph{the main term} of the
distribution $P^n$ and satisfies the local limit theorem; the
kernel  $R^n$ represents \emph{the remainder term}, and typically
decreases rapidly (see statements (iii) and (iii$'$) of Theorem
\ref{t11} below). Such kind of a representation appears to be
powerful enough to provide non-trivial applications (see Section 2
below). On the other hand, the conditions that we impose on the
distribution of $\xi_k$ in order to provide such a decomposition
to exist are very mild; in a simplest cases, these conditions have
"if and only if" form (see Theorem \ref{t13} below). Our main tool
in the current research is a certain modification of the
\emph{Malliavin calculus}.

Let us make a brief overview of the bibliography in the field.
Malliavin calculus have not been used widely for studying the properties
 of the distributions of the processes
defined by the difference relations of the type
(\ref{01}),(\ref{02}). The only paper in this direction available
to the author is \cite{Klo_Ver}, where rather restrictive
conditions are imposed both on the coefficients ($b$ should be
constant) and the distribution of $\xi_k$
 (it should possess the density from the class $C^d$). A powerful group of results is presented in the papers
\cite{kon_mammen},\cite{kon}, where a modification of
 the parametrix method in a difference set-up is developed.
  When applied to the problem formulated above, the results of \cite{kon_mammen},\cite{kon}
allow one to prove that $p^n$ converge to $p$ with the best possible rate
 $O({1\over \sqrt n})$. However, conditions imposed on the distribution of $\xi_k$ in \cite{kon_mammen},\cite{kon},
  are somewhat more restrictive than those used in our approach. For instance,
  condition (A2) of \cite{kon_mammen} requires, in our
settings, $\xi_k$ to possess the density of the class $C^4(\Re^d)$
(compare with the condition (B3) in Theorem \ref{t11} below).

The paper is organized in the following way. In Section 1, we
formulate the main theorem of the paper together with its
particular version, that is an intermediate between classic
Gnedenko's and Prokhorov's local limit theorems. In the same
section, we discuss briefly some possible improvements of the main
result. In Section 2, two applications are given. In Section 3,
the construction of the \emph{partial Malliavin calculus}, that is
our main tool, is explained in details. In Section 4, the proofs
of the main results are given.

\section{The main results}

\subsection{Formulation}
Let us introduce the notation. We write $\|\cdot\|$ for the
Euclidean norm, not indicating explicitly the space this norm is
written for. The adjoint matrix for the matrix $A$ is denoted by
$A^*$.
  The classes of functions, that have  $k$ continuous derivatives, and
functions, that are continuous and bounded together with their $k$
derivatives, are denoted by $C^k$ and $C^k_b$, correspondingly.
The derivative (the gradient) is denoted by $\nabla$, the partial
derivative w.r.t. the variable $x_r$ is denoted by $\prt_r$. The
Lebesgue measure on $\Re^d$ is denoted by $\lambda^d$. For the
measure $\mu$ on $\Bf(\Re^d),$ $\mu^{ac}$ denotes its absolutely
continuous component w.r.t. $\lambda^d$. Any time the kernel $P^n$
is decomposed into a sum $P^n=Q^n+R^n$, we mean that the kernels
$Q^n, R^n$ are non-negative; the same convention is used  for
decompositions of measures, also.
 In order to simplify notation we consider the processes
defined by (\ref{00}),(\ref{01}) and (\ref{02}) for $t\in[0,1]$
only. Of course, all the statements given below have their
straightforward analogues on an arbitrary finite time interval
$[0,T]$.

Through all the paper, $\kappa$ is a fixed integer, $\kappa\geq
4$. We denote $\epsilon(\kappa)={\kappa^2-3\kappa-2\over
2\kappa+2}$.
\begin{thm}\label{t11}  Let the following conditions hold true.

$\mathrm{(B1)}$ $a\in C^{(d+2)^2}_b(\Re^d,\Re^d), b\in
C_b^{(d+2)^2}(\Re^d,\Re^{d\times d})$ and there exists
$\beta=\beta(b)>0$ such that
$$
(b(x)b^*(x)v, v)_{\Re^d}\geq \gamma \|v\|^2, \quad x,v\in\Re^d.
$$

$\mathrm{(B2_\kappa)}$ $E\|\xi_k\|^{\kappa}<+\infty$.

$\mathrm{(B3)}$ There exist $\alpha\in (0,1)$ and bounded open set
$U\subset \Re^d$ such that
$$
{d\mu^{ac}\over d\lambda^d}\geq  {\alpha\over \lambda^d(U)}\1_U
\quad \lambda^d-\hbox{a.s.}
$$

Then $P^n$ can be represented in the form $P^n=Q^n+R^n$ in such a
way that

\begin{itemize}

\item[(i)] $Q^n_{x,t}(dy)=q_t^n(x,y)\,dy$ and $q^n\to p, n\to+\infty$ uniformly on the set
$[\delta,1]\times\Re^d\times
\Re^d$ for every $\delta\in(0,1)$;

\item[(ii)] there exist  constants  $B,C,\gamma>0$ such that, for
$t\in[0,1]$,
$$
q^n_t(x,y)\leq \begin{cases} Ct^{-{d\over
2}}\exp\brakes{-{\gamma\|x-y\|^2\over t}},& \|x-y\|\leq
tBn^{1\over \kappa+1}\\
Ct^{-{d\over 2}}\exp\brakes{-{\gamma n^{1\over \kappa+1}\|x-y\|
}},& \|x-y\|> tBn^{1\over \kappa+1}
\end{cases}
 ;
$$
in addition, for every $p>1$ there exists $C_p>0$ such that, for
$t\in[0,1], x,y\in\Re^d$,
$$
q^n_t(x,y)\leq  C_pt^{-{d\over 2}}\left(1+{\|x-y\|^2\over
t}\right)^{-p};
$$

 \item[(iii)] there exist constants $D,\rho>0$ such
that $R^n_{x,t}(\Re^d)\leq D[n^{-{\epsilon(\kappa)}}+ e^{-\rho
nt}]$, $x\in\Re^d, t\in[0,1].$
\end{itemize}

\noindent If the condition $\mathrm{(B2_\kappa)}$ is replaced by
the
 stronger condition

$\mathrm{(B2_{\exp})}$ $\exists\kap>0$ such that
$E\exp[\kap\|\xi_k\|^2]<+\infty$,

\noindent then the following stronger analogues of
$\mathrm{(ii)}$, $\mathrm{(iii)}$ hold true:
\begin{itemize}
\item[(ii$'$)] there exist  constants  $C,\gamma>0$ such that
$$
q^n_t(x,y)\leq Ct^{-{d\over 2}}\exp\brakes{-{\gamma\|x-y\|^2\over
t}},\quad t\in [0,1],\,  x,y\in \Re^d;
$$
\item[(iii$'$)] there exist  constants $D,\rho>0$ such that
$R^n_{x,t}(\Re^d)\leq De^{-\rho nt}$, $x\in\Re^d,t\in[0,1].$
\end{itemize}

\end{thm}

Let us formulate separately a particular version of Theorem
\ref{t11}. The most studied partial case of (\ref{01}) is
$a\equiv 0, b\equiv I_{\Re^d}$. In this case, $X_n\brakes{1}$ is
just the normalized sum $n^{-{1\over 2}}\sum_{k=1}^n\xi_k$ and the
limiting behavior of the distributions of such kind of a sums is
given by the CLT. For the densities of the truncated
distributions,  the following criterium can  be derived. We denote
by $P_n$ the distribution of $n^{-{1\over 2}}\sum_{k=1}^n\xi_k$.

\begin{thm}\label{t13}  The following statements are equivalent
\begin{itemize}
\item[1.] There exists $n_0\in \NN$ such that $[P_{n_0}]^{ac}$ is
not equal to zero measure.

\item[2.] There exists a representation of $P_n$ in the form
$P_n=Q_n+R_n$,  such that
\begin{itemize}
\item[(2i)] $Q_n(dy)=q_n(y)\,dy$ and
$\sup_{y\in\Re^d}\left|q_n(y)-(2\pi)^{-{d\over
2}}e^{-{\|y\|^2\over 2}}\right|\to 0, n\to \infty;$

\item[(2ii)]  there exist  constants $D,\rho>0$ such that
$R_n(\Re^d)\leq De^{-\rho n}, n\in \NN.$
\end{itemize}
\end{itemize}
\end{thm}

The well known theorem by Prokhorov states that the given above
statement 1 is equivalent to $L_1$-convergence of the density of
$[P_n]^{ac}$ to the standard normal density (see \cite{ibra_linn},
Theorem 4.4.1 for the case $d=1$). There exist examples  showing
that, even while $P_1\ll \lambda^d$, the density of $P_n$ may fail
to converge to the standard normal density uniformly (see
\cite{ibra_linn}, Ch. 4 \S3 for the example by Kolmogorov and
Gnedenko). The criterium of the uniform convergence is given by
another well known theorem by Gnedenko: such a convergence holds
if and only if there exists $n_0\in \NN$ such that $P_{n_0}$
possesses a \emph{bounded} density (see \cite{ibra_linn}, Theorem
4.3.1 for the case $d=1$).  Theorem \ref{t13} shows the following
curious feature: under condition of the Prokhorov's criterium,
some exponentially negligible remainder term can be removed from
the total distribution in such a way that, for  the truncated
distribution, the statement of the Gnedenko's theorem holds. This
feature does not seem to be essentially new; one can provide it by
using the Fourier transform technique, that is the standard tool
in the proofs of the  Prokhorov's and Gnedenko's theorems. We give
a simple proof of Theorem \ref{t13} using the partial Malliavin
calculus, developed in Section 3 below. This illustrates that the
partial Malliavin calculus is a powerful tool that allow one to
provide local limit theorems in a precise (in some cases, an "if
and only if") form.

\subsection{Some possible improvements} In the present paper, in order to keep exposition
reasonably short and transparent, we formulate the main results not in their widest possible generality.
In this subsection, we discuss shortly what kind of improvements can be made in the context of our research.

\textbf{1.} Difference relation (\ref{01}) is written w.r.t.
uniform partitions $\{0=t_n^0<t_n^1<\dots\}, t_n^k={k\over n},
k\in\ZZ^+,n\in \NN$. Without a significant change of the proofs,
one can prove analogues of Theorem \ref{t11} for the processes,
defined by the difference relations of the type (\ref{01}) with
${\xi_k\over \sqrt{n}}$ replaced by
$\xi_k\cdot\sqrt{t_n^k-t_n^{k-1}}$ and partitions $\{t_n^k\}$
satisfying condition \be\label{11} \exists\, c,C,d,D>0:\quad
\lim\inf_{n\to+\infty}{1\over n}\#\left\{k\Big|t_n^k\leq
t,(t_n^k-t_n^{k-1})\in \left[{c\over n},{C\over
n}\right]\right\}\geq td,
$$
$$
 \lim\sup_{n\to+\infty}{1\over
n}\#\left\{k\Big|t_n^k\leq t,(t_n^k-t_n^{k-1})\in \left[{c\over
n},{C\over n}\right]\right\}\leq tD,\quad t\in(0,1].
 \ee

\textbf{2.} One can, without a significant change of the proofs,
replace  the sequence of i.i.d. random vectors $\{\xi_k\}$ in
(\ref{01}) by a triangular array $\{\xi_{n,k}, k\leq n\}$ of
independent random vectors, possibly not identically distributed,
having zero mean and identity for the covariance matrix. Under
such a modification, condition (B2$_\kappa$) should be replaced by
$\sup_{n,k}E\|\xi_{n,k}\|^\kappa<+\infty$, and condition (B3)  by
$$
\exists\, \alpha,r>0, x_n\in \Re^d:\quad {d[\mu_{n,k}]^{ac}\over
d\lambda^d }\geq \alpha \1_{B(x_n,r)}\quad \lambda^d-\hbox{a.s.},
 \leqno\mathrm{(B3')}
$$
here $\mu_{n,k}$ denotes the distribution of $\xi_{n,k}$, $B(x,r)$
denotes the open ball in $\Re^d$ with the  centrum $x$ and radius
$r$. Also, the phase space for $\xi_{n,k}$ may be equal $\Re^m$
with $m\geq d$ (note that the case $m<d$ is excluded by the
condition (B1)).

\textbf{3.} Under an appropriate regularity conditions on $a,b$,
Malliavin's representation, analogous to (\ref{315}), can be
written for the derivatives of the truncated density of an
arbitrary order with respect to both $x$ and $y$. Thus, after some
standard technical steps, one can obtain the following estimate,
that generalize statement (ii$'$) of Theorem \ref{t11}: for a
given $k,l\in \NN$,
$$
{\prt^{k+l}\over \prt {x}^k\prt y^l}q_t^n(x,y)\leq
C_{k+l}t^{-{d+k+l\over 2}}\exp\left(-{\gamma\|y-x\|^2\over
t}\right)
$$
under (B1), (B2$_{\exp}$),(B3) and $a\in
C^{(d+k+l+1)^2}_b(\Re^d,\Re^d),$ $b\in
C_b^{(d+k+l+1)^2}(\Re^d,\Re^{d\times d})$.

\textbf{4.}
 Theorem \ref{t31} provides the truncated limit theorem
without essential restrictions on the structure of the
functionals. For instance, one can apply this theorem in order to
obtain a truncated local limit theorem for difference
approximations of integral functionals, etc.

\textbf{5.} Like the Malliavin calculus for (continuous time)
diffusion processes, the partial Malliavin calcucus, developed in
Section 3,  can be applied when the diffusion matrix is not
uniformly elliptic, but  locally elliptic, only. However, the
changes that should be done in the proof are significant; in
particular, Theorem \ref{t31} is not powerful enough to cover this
case. Thus we postpone the detailed investigation of this case
(and more generally, the case of coefficients satisfying an
analogue of H\"ormander condition) to some further research.

\textbf{6.} In the present paper, we concentrate on the individual
estimates for the densities $q_t^n$ and do not deal with the
convergence rate in the statement (i) of Theorem \ref{t11}.  The
(seemingly) possible way to establish such a rate is to write the
Malliavin's representation, analogous to (\ref{315}), for the
limiting density $p$ and then construct both the functionals
$X_n(t)\equiv f_n$ and $X(t)=f$ and the corresponding
\emph{weights} $\Upsilon^{f_n}$ and $\Upsilon^f$, involved into
the Malliavin's representation, on the same probability space with
a controlled $L_2$-distance between $(f_n,\Upsilon^{f_n})$ and
$(f,\Upsilon^f)$. Since the question about the estimates in the
\emph{strong invariance principle} for the pair
$(f_n,\Upsilon^{f_n})$ is far from being trivial, we postpone the
detailed investigation of the rate of convergence in Theorem
\ref{t11} to some further research. We remark that the
modification of
 the parametrix method, developed in
\cite{kon_mammen},\cite{kon}, provides, under more restrictive
conditions on the distribution of $\{\xi_k\}$, the best possible
converence rate $O({1\over \sqrt n})$.

\section{Applications}

In this section, we formulate two applications of Theorem
\ref{t11}.

\subsection{Mixing and
convergence rates for difference approximations to SDE's}

 Under condition (B1) and some recurrence conditions,
the process $X$ is \emph{ergodic}, i.e.  possesses a unique
invariant  distribution $\mu_{inv}$ (see
\cite{veretennikov_1},\cite{veretennikov_2}). Moreover, an
explicit estimates for the $\beta$-\emph{mixing coefficients} and
for the rate of convergence of $P_{x,t}\equiv P(X(t)\in
\cdot|X(0)=x)$ to $\mu_{inv}$ in total variation norm are also
available.

The processes $X_n$, restricted to ${1\over n}\ZZ_+$, are a Markov
chains. The following natural question takes its origins in a
numerical applications: can the mentioned above estimates for the
mixing and convergence rate be made uniform over the class
$\{X_n,n\geq 1, X\}$? This question is studied in the recent paper
\cite{Klo_Ver}, see more discussion therein.

 In this subsection, we use the truncated local limit theorem (Theorem \ref{t11})
in order to establish the required uniform estimates. Denote, by
$\|\cdot\|_{var}$, the total variation norm. Recall that the
$\beta$-mixing coefficient for $X$ is defined by
$$
\beta_{x}(t)\equiv \sup_{s\in\ax}\Es_\mu\left\|P(\cdot|\Ff_0^s,
X(0)=x)-P(\cdot|X(0)=x) \right\|_{var,\Ff_{t+s}^\infty},\quad
t\in\ax,
$$
where $\Ff_a^b\equiv \sigma (X(s),s\in [a,b]), P(\cdot|\Ff_0^s,
X(0)=0)$ denotes the conditional distribution of the process $X$
with $X(0)=x$ w.r.t. $\Ff_0^s$, and
$$
\left\|\varkappa\right\|_{var,\Gf}\eqdef \sup_{\begin{array}{c}
B_1\cap B_2=\emptyset,B_1\cup B_2=\Re^m\\
B_1,B_2\in\Gf
\end{array}}[\varkappa(B_1)-\varkappa(B_2)].
$$
The $\beta$-mixing coefficient $\{\beta^n_\mu(t), t\in {1\over
n}\ZZ_+\}$ for $X_n$ is defined analogously.

\begin{thm}\label{t22}

Let conditions $\mathrm{(B1)}$ and $\mathrm{(B3)}$ hold true.
Suppose also  that
\begin{itemize}
\item[(B4)] there exists $R_0>0$ and $r>0$ such that
$$
(a(x),x)_{\Re^d}\leq -r\|x\|, \quad \|x\|\geq R_0.
$$
\item[(B5)] there exists $\kap>0: E\exp[\kap\|\xi\|]<+\infty.$

\end{itemize}

Then, for every process $X, X_n,n\geq 1$ there exists unique
invariant distribution $\mu_{inv}, \mu_{inv}^n$. Moreover, there
exist $n_0\in \NN$, a function $C:\Re^d\to \ax$ and a constant
$c>0$ such that
$$
\|P_{x,t}^n-\mu_{inv}^n\|_{var}\leq C(x)e^{-ct},\quad t\in {1\over
n}\ZZ_+,n\geq n_0, \quad  \|P_{x,t}-\mu_{inv}\|_{var}\leq
C(x)e^{-ct},\quad t\in\ax,
$$
$$
\beta_x^n(t)\leq C(x)e^{-ct},\quad t\in {1\over n}\ZZ_+,n\geq
n_0,\quad \beta_x(t)\leq C(x)e^{-ct},\quad t\in\ax.
$$
\end{thm}

{\it Remarks. 1.} The statement of  Theorem \ref{t22} is analogous
to the one of Theorem 1 \cite{Klo_Ver}. The main improvement is
that the conditions ($D_1$) -- ($D_3$) of Theorem 1 \cite{Klo_Ver}
are replaced by (seemingly, the mildest possible) condition (B3).
In addition, Theorem \ref{t22}, unlike Theorem 1 \cite{Klo_Ver},
admits non-constant diffusion coefficients $b$.

{\it 2.} The mixing and convergence rates established in Theorem
\ref{t22} are called an \emph{exponential} ones. If the recurrence
condition (B4) is replaced by a weaker ones, then the
\emph{subexponential} or \emph{polynomial} rates can be
established (see Theorem 1 \cite{Klo_Ver}, cases 2 and 3). We do
not give an explicit formulation here in order to shorten the
exposition.

\subsection{Difference approximation for local times of multidimensional
diffusions}

Consider a \emph{W-measure} $\mu$ on $\Re^d$, that is, by
definition (\cite{dynkin}, Chapter 8), a $\sigma$-finite measure
satisfying the condition \be\label{51}
\sup_{x\in\Re^m}\int_{\|y-x\|\leq
1}w_d(\|y-x\|)\mu(dy)<+\infty \quad \hbox{ with }\quad  w_d(r)= \begin{cases} r, &d=1\\\max(-\ln r, 1),& d=2\\
r^{2-d},& d>2
\end{cases}.\ee

Every such a measure generates a \emph{W-functional}
(\cite{dynkin}, Chapter 6) of a Wiener process $W$ on $\Re^d$,
\be\label{52} \phi^{s,t}=\phi^{s,t}(W)=\int_s^t{d\mu\over
d\lambda^d}(W(r))\, dr, \quad 0\leq s\leq t. \ee  For singular
$\mu$, equality (\ref{52}) is a formal notation, that can be
substantiated via an approximative procedure with $\mu$
approximated by an absolutely continuous measures (\cite{dynkin},
Chapter 8). The functional $\phi$ is naturally interpreted as the
\emph{local time} for the Wiener process, correspondent to the
measure $\mu$.

Next, let the process $X$ be defined by (\ref{00}) and satisfy
(\ref{A2}), that means that the asymptotic behavior of its
transition probability density as $t\to 0+$ is similar to the one
of the transition probability density for the Wiener process. Then
the estimates, analogous to those given in \cite{dynkin}, Chapter
8 provide that the $W$-functional of the process $X$ \be\label{53}
\phi^{s,t}=\phi^{s,t}(X)=\int_s^t{d\mu\over d\lambda^d}(X(r))\,
dr, \quad 0\leq s\leq t \ee is well defined. We interpret this
functional as the \emph{local time} for the diffusion process $X$,
correspondent to the measure $\mu$.

At last, let the sequence $X_n, n\in \NN$ of difference
approximations for $X$ be defined by (\ref{01}),(\ref{02}).
 Consider a sequence of the
functionals $\phi_n(X_n)$ of the processes $X_n$ of the form
\be\label{54}\phi_{n}^{s,t}=\phi_{n}^{s,t}(X_n)\eqdef {1\over
n}\sum_{k:s\leq {k\over
n}<t}F_{n}\brakes{X_n\brakes{\frac{k}{n}}},\quad 0\leq s<t.\ee In
Theorem \ref{t21} below, we establish  sufficient conditions
 for the joint distributions of
$(\phi_n,X_n)$  to converge weakly to the joint distribution of
$(\phi,X)$. Thus, it is natural to say that the functionals
$\phi_n$ defined by (\ref{54}) provide the difference
approximation for the local time $\phi$ defined by (\ref{53}). For
the further discussion and references concerning this problem, we
refer the reader to the recent paper \cite{kulik_da}.

We fix $x\in\Re^d$ and suppose that $X_n(0)=X(0)=x$. We denote
$\TT=\{(s,t):0\leq s\leq t\}$. In order to shorten exposition, we
suppose $\mu$ to be finite and to have a compact support. Together
with the functionals $\phi_n$ that are discontinuous w.r.t.
variables $s,t$, we consider the "random broken line" processes
$$
\psi_n^{s,t}=\phi_{n}^{{j-1\over n}, {k-1\over
n}}-(ns-j+1)\phi_{n}^{{j-1\over n}, {j\over
n}}+(nt-k+1)\phi_{n}^{{k-1\over n}, {k\over n}},\quad
s\in\left[{j-1\over n},{j\over n}\right),t\in\left[{k-1\over
n},{k\over n}\right).
$$

\begin{thm}\label{t21} Let conditions $\mathrm{(B1)}$, $\mathrm{(B2_{6})}$, $\mathrm{(B3)}$ hold true.
Suppose also  that
\begin{itemize}
\item[(B6)] $F_n(x)\geq 0, x\in\Re^d, n\geq 1$ and ${1\over
n}\sup\limits_{x\in\Re^d}F_n(x)\to 0, n\to\infty$;

\item[(B7)] Measures $\mu_n(dx)\equiv F_n(x)\lambda^d(dx)$ weakly
converge to $\mu$;
\end{itemize}
$$
\lim_{\delta\downarrow
0}\mathop{\lim\sup}_{n\to+\infty}\sup_{x\in\Re^d}\int_{\|y-x\|\leq
\delta}w_d(\|y-x\|)\mu_n(dy)\to 0. \leqno\quad \mathrm{(B8)}
$$

Then $(X_n,\psi_n(X_n))\Rightarrow (X,\phi(X))$ in a sense of weak
convergence in $C(\ax, \Re^d)\times C(\TT,\Re^+)$.
\end{thm}

{\it Remarks. 1.} The statement of  Theorem \ref{t21} is analogous
to the one of Theorem 2.1 \cite{kulik_da}. The main improvement is
that the condition A3) of Theorem 2.1 \cite{kulik_da} is replaced
by (seemingly, the mildest possible) condition (B3).

{\it 2.} Once Theorem \ref{t21} is proved, one can use the
standard truncation procedure in order to replace the moment
condition $\mathrm{(B2_{6})}$ by the Lyapunov type condition
"$\exists\,\delta>0: E\|\xi_k\|^{2+\delta}<+\infty$" (e.g.
\cite{kar_kul} Section 5).

{\it 3.} For  examples and a discussion on the relation between
conditions (\ref{51}) and (B8), we refer the reader to
\cite{kulik_da}.

\section{Partial Malliavin calculus on a space with a product measure}

 For every given $n\in \NN$ and $t\in[0,1]$, the value $X_n(t)$ is a functional of
$\xi_1,\dots, \xi_n$ and thus can be interpreted as a functional
on the space $(\Re^d)^n$ with the product measure $\mu^n$.
However, under the conditions of Theorem \ref{t11},  $\mu$ may
contain a singular component and therefore it may fail to have
logarithmic derivative. Thus, in general, one can not write the
integration-by-parts formula on the probability space
$((\Re^d)^n,(\Bf(\Re^d))^{\otimes n}, \mu^n)$. We overcome this
difficulty by using the following trick. Under condition (B3), the
measure $\mu$ can be decomposed into a sum \be\label{31}
\mu=\alpha \cdot \pi_U+(1-\alpha)\cdot \nu, \ee where $\pi_U$ is
the uniform distribution on $U$.  One can write (on an appropriate
probability space) the representation for $\{\xi_k\}$
corresponding to (\ref{31}): \be\label{32} \xi_k=\eps_k\cdot
\eta_k+(1-\eps_k)\cdot \zeta_k, \ee where $\eta_k\sim \pi_U,
\zeta_k\sim \nu,$ and   the distribution $\varkappa$ of $\eps_k$
is equal to Bernoulli distribution with $\varkappa\{1\}=\alpha$.
This representation allows one to consider the family
$\xi_1,\dots, \xi_n$ (and, therefore, the process $X_n$) as a
functional on the following probability space: \be\label{3200}
\Omega=(\Re^d\times\{0,1\}\times\Re^d)^n, \quad
\Ff=(\Bf(\Re^d)\otimes 2^{\{0,1\}}\otimes\Bf(\Re^d))^n, \quad
P=(\pi_U\times\kap\times \nu)^n.\ee Now, the measure $\pi_U$ has a
logarithmic derivative w.r.t. a properly  chosen vector field, and
some kind of an integration-by-parts formula can be written on the
probability space $(\Omega,\Ff,P)$ (see subsection 3.1 below). The
Malliavin-type calculus, associated to this formula, is our main
tool in the proof of Theorems \ref{t11}, \ref{t13}. We call this
calculus a \emph{partial} one because the stochastic derivative,
this calculus is based on, is defined w.r.t. a proper group of
variables, while the other variables play the role of interfering
terms. In this section, we give the main constructions of the
partial Malliavin calculus, associated to the representation
(\ref{32}).

\subsection{Integration-by-parts formula. Derivative and divergence. Sobolev classes.}

Denote $\Omega=\Omega_1\times\Omega_2\times\Omega_3$,$$
 \Omega_1=\Omega_3=(\Re^d)^n, \quad
\Omega_2=\{0,1\}^n. $$ We write a point $\omega\in \Omega$  in the
form $\omega=(\eta,\eps,\zeta)$, where
$$
\eta=(\eta_1,\dots, \eta_n)\in (\Re^d)^n,\quad
\eps=(\eps_1,\dots,\eps_n)\in \{0,1\}^n,\quad
\zeta=(\zeta_1,\dots, \zeta_n)\in (\Re^d)^n
$$
and $\eta_k=(\eta_{k1},\dots,
\eta_{kd}),\zeta_k=(\zeta_{k1},\dots, \zeta_{kd})$. In this
notation, the random variables $\eta_k,\eps_k,\zeta_k$ are defined
just as the coordinate functionals:
$$
\eta_k(\omega)=\eta_k,\,\eps_k(\omega)=\eps_k,\,\zeta_k(\omega)=\zeta_k,\quad
\omega=(\eta,\eps,\zeta)\in\Omega.
$$
Denote by $\Cf$ the set of bounded measurable functions $f$ on
$\Omega$ such that, for every $(\eps,\zeta)\in \Omega_2\times
\Omega_3$, the function $f(\cdot,\eps,\zeta)$ belongs to the class
$C^\infty(\Re^d)$ and
$$
\mathop{\mathrm{ess}\sup}_{\eta,\eps,\zeta} \|[\nabla_\eta]^j
f((\eta,\eps,\zeta))\|<+\infty, \quad j\in \NN,
$$
where $\nabla_\eta$ denotes the gradient w.r.t. variable $\eta$.

For $f\in \Cf$ and $k=1,\dots, n, r=1,\dots, d$, denote by
$\prt_{kr} f$ the derivative of $f$ w.r.t. the variable
$\eta_{kr}$. Also, denote by $H$ the space $\Re^{d\times n}$
considered as a (finite-dimensional) Hilbert space with the usual
Euclid norm, and by $\{e_{kr}, r=1,\dots, d, k=1,\dots,n\}$ the
canonical basis in it: all coordinates of the vector $e_{kr}$ are
equal to zero except the coordinate with the index $kr$ being
equal to one. For a given functions $\psi:\Re^d\to \Re$ and
$\theta_n:(\Re^d)^n\to [0,1]$, define the stochastic gradient $D$
by the formula \be\label{320} [D
f](\eta,\eps,\zeta)=\theta_n(\zeta)\sum_{k,r}\psi(\eta_k)\cdot
[\prt_{kr}f](\eta,\eps,\zeta)\cdot e_{kr},\quad f\in \Cf. \ee This
definition can be naturally extended to the functionals taking
their values in a finite-dimensional  Hilbert space $Y$ (actually,
in any separable Hilbert space, but we do not need such a
generality in our further construction). Given an orthonormal
basis $\{y_l\}$ in $Y$, denote by $\Cf^Y$ the set of the functions
of the type
$$
y=\sum_lf_l\cdot y_l,\quad \{f_l\}\subset \Cf
$$
and put for such a function
$$
Dy=\sum_l[Df_l]\otimes y_l.
$$
It is easy to see that the definitions of the class $\Cf^Y$ and
the derivative $D$ do not depend on the choice of the basis
$\{y_l\}$. By the construction, $D$ satisfies the chain rule: for
any two spaces $Y, Z$ and for any  $f_1,\dots,f_m\in \Cf^Y,\, F\in
C^\infty(Y^m, Z),\, m\geq 1$,  \be\label{321} F(f_1,\dots,f_m)\in
\Cf^Z\quad \hbox{and}\quad
D\Big[F(f_1,\dots,f_m)\Big]=\sum_{j=1}^m [\prt_j
F](f_1,\dots,f_m)\cdot Df_j.
 \ee
We denote $D^0f=f, D^1f=Df$. The higher derivatives $D^j, j>1$ are
defined iteratively: $D^j=\mathop{\underbrace{D\cdot\dots \cdot
D}}\limits_{j}$ (note that the first operator in this product acts
on the elements of $\Cf^Y$ while the last one acts on the elements
of $\Cf^{H^{\otimes (j-1)}\otimes Y}$).

Everywhere below, we suppose that $U$ is an open ball $B(z,r)$
(this obviously does not restrict generality). We define the
function $\psi$ in (\ref{320}) by $\psi(x)=r^2-\|x-z\|^2$. Due to
this choice, $\psi\in C^\infty(\Re^d)$ and $\psi=0$ on $\prt U$.
These properties of $\psi$ imply the following integration by
parts formula:
$$
\int_{U}[\prt_r f](x) \psi(x)\, dx=-\int_U [\prt_r \psi](x)f(x)\,
dx,\quad f\in C^1(\Re^d), \,r=1,\dots,d.
$$
As a corollary of this formula, we obtain the following statement.

\begin{prop}\label{p31} For every $h\in H$ and every $f\in \Cf$, the following
integration-by-parts formula holds true:
  \be\label{33} E (Df,h)_H =-E (\rho,h)_H f,\quad  \rho\equiv
\theta_n(\zeta)\cdot\sum_{k,r}[\prt_r\psi](\eta_k)\cdot
  e_{kr}.
  \ee
  \end{prop}

The formula (\ref{33}) allows one to introduce, in a standard way,
the divergence operator corresponding to the derivative $D$. For
$g\in \Cf^{H\otimes Y}$, put \be\label{34}
\delta(g)=-\sum_{k,r,l}\Big[(\rho,e_{kr}) g_{krl}+(D g_{krl},
e_{kr})_H\Big]\cdot y_l, \quad g_{krl}=(g,e_{kr}\otimes
y_l)_{H\otimes Y}. \ee By the choice of the function $\psi$,
 $\delta(g)\in \Cf^Y$ as soon as $g\in \Cf^{H\otimes
Y}$. The chain rule (\ref{321}) and the integration-by-parts
formula (\ref{33}) imply that the operators $D$ and $\delta$ are
mutually adjoint in a sense of the following \emph{duality
formula}: \be\label{35} E(Df, g)_H=Ef\delta(g),\quad f\in \Cf^Y,
\,g\in\Cf^{Y,H}.\ee Since, for every $p\geq 1$ and every $Y$,
$\Cf^Y$ is dense in $L_p(\Omega,P,Y)$, the duality formula
(\ref{35}) provides that, for any $p\geq 1$ and any $Y$, the
operators $D,\delta$ are closable as densely defined unbounded
operators
$$
D:L_p(\Omega,P,Y)\to L_p(\Omega,P,H\otimes Y), \quad \delta:
L_p(\Omega,P,H\otimes Y)\to L_p(\Omega,P,Y).
$$

\begin{dfn}\label{d31} The Sobolev class $W^m_p(Y)$ ($p\geq 1, m\in \ZZ_+$) is the completion of
the class $\Cf^Y$ w.r.t. the norm
$$
\|f\|_{p,m}\equiv\left[\sum_{j=0}^m E\|D^j f\|_{ H^{\otimes
j}\otimes Y}^p\right]^{1\over p}<+\infty.
$$
Since $D$ is closable in $L_p$ sense, there exists the  canonical
embedding of $W_p^m(Y)$ into $L_p(\Omega,P,Y)$.
\end{dfn}
We denote $W_\infty^\infty(Y)=\bigcap_{m,p}W_p^m(Y)$. If $Y=\Re$
then we denote the corresponding Sobolev spaces simply by $W_p^m$.

\subsection{Algebraic relations for derivative and divergence. Moment estimates.}

Let us introduce some notation. We denote by $\Cs$ a constant such
that its value can be calculated explicitly, but this calculation
is omitted. The value of $\Cs$ may vary from line to line. If the
value of the constant $\Cs$ depends on some parameters, say $m,d$,
then we write $\Cs(m,d)$.  The latter notation indicates that the
value of the constant does not depend on other parameters (for
instance, $n$). If, in a sequel, the constant $\Cs$ is referred
to, then we endow it with the lower index like $\Cs_1,\Cs_2$, etc.
We use standard notation $\{\delta_{jk},j,k\in\NN\}$ for the
Kroeneker's symbol.

  For an
$H\otimes H\otimes Y$-valued element $K$, we denote by $K^*$ the
element such that
$$
(K^*,h\otimes g\otimes y)_{H\otimes H\otimes Y}=(K,g\otimes
h\otimes y)_{H\otimes H\otimes Y}, \quad h,g\in H,y\in Y.
$$
For an $X\otimes Y$-valued element $g_1$ and $X\otimes Z$-valued
element $g_2$, we denote by $(g_1,g_2)_{X}$ the $Y\otimes
Z$-valued element
$$
(g_1,g_2)_{X}\equiv \sum_{l_1,l_2,l_3}(g_1,x_{l_1}\otimes
y_{l_2})_{X\otimes Y}\cdot (g_2,x_{l_1}\otimes z_{l_3})_{X\otimes
Z}\cdot [y_{l_2}\otimes z_{l_3}],
$$
here $\{x_l\},\{y_l\}$, $\{z_l\}$ are orthonormal bases in $X,Y$
and $Z$, correspondingly. We also denote for an $Y\otimes
X$-valued element $g_1$ and $Z\otimes X$-valued element $g_2$
$$
(g_1,g_2)_{X}\equiv \sum_{l_1,l_2,l_3}(g_1,y_{l_2}\otimes
x_{l_1})_{Y\otimes X}\cdot (g_2,z_{l_3}\otimes x_{l_1})_{Z\otimes
X}\cdot [y_{l_2}\otimes z_{l_3}].
$$
Although the same notation $(\cdot,\cdot)_X$ is used for two
slightly different objects,  it does not cause misunderstanding
further.

Consider  the $\Lf(H)$-valued random element (i.e., \emph{random
operator} in $H$) $B$, defined by the relations \be\label{351}
\Big(B e_{k_1 r_1},e_{k_2r_2}\Big)_H = -\Big(D (\rho,e_{k_1
r_1})_H,e_{k_2r_2}\Big)_H=-{\theta_n^2(\zeta)}\cdot \delta_{k_1
k_2}\cdot [\prt_{r_1}\prt_{r_2}\psi](\eta_{k_1})\cdot \psi
(\eta_{k_1}), \ee $ k_{1,2}=1,\dots, n, \, r_{1,2}=1,\dots, d$. We
define the action of $B$ on $H\otimes Y$-valued element $g$ by
$$
Bg=\sum_{k,r,l} (g,e_{k,r}\otimes y_l)_{H\otimes Y}\cdot\Big[
[Be_{k,r}]\otimes y_l\Big].
$$

 Using representation (\ref{34}), one can deduce the
commutation relations for the operators $D,\delta$, analogous to
those for the stochastic derivative and integral for the Wiener
process (the proof is straightforward and omitted; for the Wiener
case, see \cite{nualart}, \S 1.2).
\begin{prop}\label{p32} $\mathrm{I.}$ If $f\in \Cf, g\in \Cf^{H\otimes
Y}$, then $f\cdot g \in \Cf^{H\otimes Y}$ and $$ \delta(f\cdot
g)=f\cdot \delta (g)-(Df, g)_H.$$

$\mathrm{II.}$ If $g\in \Cf^{H\otimes Y}$, then
$$
D\Big[\delta(g)\Big]=Bg+\delta\Big([Dg]^*\Big).
$$

$\mathrm{III.}$ If $g_1,g_2\in \Cf^H$, then
$$
(D[\delta(g_1)],g_2)_H=(Bg_1,g_2)_H+\delta\Big((Dg_1,
g_2)_H\Big)+([Dg_1]^*, Dg_2)_{H\otimes H}.
$$
\end{prop}

The main result of this subsection is given by the following
lemma.

\begin{lem}\label{l31} Let $m,l\in \NN, g\in W_{2m}^{2m+l-1}(H)$. Then there exists
$\delta(g)\in W_{2m}^l$ and \be\label{36} \|\delta(g)\|_{2m,l}\leq
\Cs(m,l,d,\psi)\cdot \|g\|_{2m,2m+l-1}. \ee
\end{lem}

\emph{Remarks. 1.} On the Wiener space, the typical way to prove
estimates of the type (\ref{36}) is to use Meyer's inequalities
for the generator $L=\delta\cdot D$ of the Ornstein-Uhlenbeck
semigroup (see, for instance, \cite{nualart} \S 2.4). Moreover, on
the Wiener space, (\ref{36}) can be made more precise: the similar
inequality holds with $2m+l-1$ replaced by $l+1$. In our settings,
it is not clear whether the operator $\delta\cdot D$ provides the
analogues of Meyer's inequalities, since it does not have the
specific structural properties of the Ornstein-Uhlenbeck generator
(such as Mehler's formula, hypercontractivity of the associated
semigroup, etc.). Thus we prove (\ref{36}) straightforwardly by
using an iterative integration-by-parts procedure.

\emph{2.} Throughout the exposition, the function $\psi$ is fixed
together with the set $U=B(x,z)$. However, when the constant $\Cs$
depends on the values of $\psi$ or its derivatives, we indicate it
explicitly  in the notation for $\Cs$.

In order to prove Lemma \ref{l31}, we need some auxiliary
statements and notation. For $g\in \Cf^Y$ and $m\in \ZZ_+$, we
define the random variable $|g|_m$ by
$$
|g|_m\equiv \left(\sum_{j=0}^m\|D^j g\|^2_{H^{\otimes j}\otimes
Y}\right)^{1\over 2}.
$$

\begin{lem}\label{l32} If $g\in \Cf^{H\otimes Y}$ then $Bg\in \Cf^{H\otimes Y}$ and,
for every $m\in \ZZ_+$,
$$
|Bg|_m\leq \Cs(m,d,\psi)\cdot |g|_{m}.
$$
\end{lem}

{\it Proof.} Write $Bg$ in the coordinate form:
$$
Bg=\sum_{k,r^1,r^2,l} (g,e_{kr^1}\otimes y_l)_{H\otimes Y}\cdot
b_{k,r^1,r^2}\cdot [e_{kr^2}\otimes y_l],
$$
where $b_{k, r^1,r^2}=(Be_{kr^1}, e_{kr^2})_H$ (recall that
$(Be_{k^1r^1}, e_{k^2r^2})_H=0$ as soon as $k^1\not=k^2$). Write
the Leibnitz formula for the higher derivatives:
$$ D^m(Bg)=\sum_{k,r^1,r^2,l}\sum_{\Theta\in
2^{\{1,\dots,m\}}}\sum_{k_1,r_1,\dots,k_m,r_m}\Big(D^{\#\Theta}(g,e_{kr^1}\otimes
y_l)_{H\otimes Y}, \bigotimes_{i\in
\Theta}e_{k_ir_i}\Big)_{H^{\otimes \#\Theta }}\times
$$
\be\label{37}
 \times \Big(D^{m-\#\Theta}b_{kr^1r^2}, \bigotimes_{i\not\in
\Theta}e_{k_ir_i}\Big)_{H^{\otimes
(m-\#\Theta)}}\Big[\bigotimes_{i=1}^m e_{k_ir_i}\otimes
e_{kr^2}\otimes y_l\Big], \ee where $\#\Theta_m$ denotes the
number of elements in the set $\Theta$. We write
$$
S_\Theta=
\sum_{k,r^1,r^2,l,k_1,r_1,\dots,k_m,r_m}\Big(D^{\#\Theta}(g,e_{kr^1}\otimes
y_l)_{H\otimes Y}, \bigotimes_{i\in
\Theta}e_{k_ir_i}\Big)_{H^{\otimes \#\Theta }}\times
$$
$$
\times \Big(D^{m-\#\Theta}b_{kr^1r^2}, \bigotimes_{i\not\in
\Theta}e_{k_ir_i}\Big)_{H^{\otimes
(m-\#\Theta)}}\Big[\bigotimes_{i=1}^m e_{k_ir_i}\otimes
e_{kr^2}\otimes y_l\Big]
$$
and estimate $\|S_\Theta\|_{H^{\otimes(m+1)}\otimes Y}$. The
function $\psi$ belongs to $C^\infty$ and is bounded together with
all its derivatives on $U$. Thus, one can deduce from
 the representation (\ref{351}) and formula (\ref{320})  that
$$
\|D^{M}b_{kr^1r^2}\|_{H^{\otimes M}}\leq \Cs(M,d,\psi),\quad M\in
\NN.
$$
In addition, due to (\ref{351}),
$$
 \Big(D^{m-\#\Theta}b_{kr^1r^2},
\bigotimes_{i\not\in \Theta}e_{k_ir_i}\Big)_{H^{\otimes
(m-\#\Theta)}}=0
$$
as soon as $k_i\not=k$ for some $i\not\in \Theta$. Using these
facts, we deduce that
$$
\sum_{k_i\in \{1,\dots, n\}, r_i\in \{1,\dots, d\}, i\not\in
\Theta} \Big(D^{m-\#\Theta}b_{kr^1r^2}, \bigotimes_{i\not\in
\Theta}e_{k_ir_i}\Big)^2_{H^{\otimes(m-\#\Theta)}}\leq
\Cs(m,d,\psi).
$$
Thus
$$
\|S_\Theta\|^2_{H^{\otimes(m+1)}\otimes Y}\leq
$$
$$\leq  \Cs(m,d,\psi)
\sum_{k,r^1,r^2,l,k_i\in \{1,\dots, n\}, r_i\in \{1,\dots, d\},
i\in \Theta}\Big(D^{\#\Theta}(g,e_{kr^1}\otimes y_l)_{H\otimes Y},
\bigotimes_{i\in \Theta}e_{k_ir_i}\Big)^2_{H^{\otimes
\#\Theta}}\leq
$$
$$
\leq \Cs(m,d,\psi)|g|_{\#\Theta}^2\leq \Cs(m,d,\psi)|g|_{m}^2.
$$
Taking the sum over $\Theta\in 2^{\{1,\dots,m\}}$ and using the
Cauchy inequality, we obtain the required statement. The lemma is
proved.

\begin{prop}\label{p33} $\mathrm{I.}$ Let  $ g_1\in \Cf^{X\otimes Y}, g_2\in
\Cf^{X\otimes Z}$, then $(g_1,g_2)_{X}\in \Cf^{Y\otimes Z}$ and
$$
|(g_1,g_2)_{X}|_m\leq \Cs(m)|g_1|_{m}|g_2|_m, \quad m\geq 0.
$$
$\mathrm{II.}$ Let $g\in \Cf^Y$ and $A\in \Lf(Y,Z)$,  then $Ag\in
\Cf^Z$ and $|Ag|_m\leq \|A\|\cdot |g|_m, m\geq 0$.
\end{prop}

The second statement is a straightforward corollary of the chain
rule (\ref{321}). The first one can be proved using the Leibnitz
formula; the proof is totally analogous to the one of Lemma
\ref{l32}, and thus we omit the detailed exposition.

\emph{Remark.} Taking $X=\Re, $ we obtain that, for $ g_1\in
\Cf^{Y}, g_2\in \Cf^{Z}$, $g_1\otimes g_2\in \Cf^{Y\otimes Z}$
with $|g_1\otimes g_2|_m\leq \Cs(m)|g_1|_{m}|g_2|_m, \quad m\geq
0.$

Using iteratively statement II of Proposition \ref{p32}, we obtain
that, for $g\in \Cf^{H}$ and $m\geq 1$, the derivative
$D^m[\delta(g)]$ can be expressed in the form
$$
D^m[\delta(g)]=F_m(g)+\delta(G_m(g)),
$$
where $F_m(g)\in \Cf^{H^{\otimes m}}, G_m(g)\in \Cf^{H^{\otimes
(m+1)}}$ are defined via the iterative procedure
$$
F_0(g)=0,\quad G_0(g)=g,\quad G_{i+1}(g)=[DG_i(g)]^*,\quad
F_{i+1}(g)=DF_i(g)+BG_i(g),\quad i\geq 0.
$$
The mapping $K\mapsto K^*$ is an isometry in $H\otimes H\otimes
Y$, thus statement II of Proposition \ref{p33} provides that
\be\label{370} |G_m(g)|_{j}\leq |g|_{m+j},\quad m,j\geq 0.\ee
Using Lemma \ref{l32}, we deduce that \be\label{371}
|F_m(g)|_j\leq \Cs(m,j,d,\psi)\cdot |g|_{m+j}, \quad m,j\geq 0.
\ee

Thus, in order to prove inequality  (\ref{36}) for $g\in \Cf^H$,
it is sufficient to prove that \be\label{372}
E\|\delta(g)\|_Y^{2m}\leq \Cs(m,d,\psi) E |g|_{2m-1}^{2m}\ee for
any $m\geq 1, g\in \Cf^{H\otimes Y}$ and arbitrary Hilbert space
$Y$. In order to prove estimate (\ref{372}) we embed it into a
larger family of estimates. Consider the following objects.

1. Numbers $k_0,\dots, k_v\in \ZZ_+$ ($v\leq 2m$) such that
$k_0+\dots+k_v=2m$. Denote
$I_{j}=[k_0+\dots+k_{j-1}+1,k_0+\dots+k_j]\cap \NN, j=1,\dots, v$,
$I_0=[1,k_0]\cap \NN$ (if $k_0=0$ then $I_0=\emptyset$).

2. Function $\sigma:\{1,\dots,2m\}\to \{1,\dots,m\}$ such that
$\#\sigma^{-1}(\{i\})=2$ and $\#\Big[I_j\cap
\sigma^{-1}(\{i\})\Big]\leq 1$ for every $i=1,\dots,m$ and $j=0,
\dots,v$.

\begin{lem}\label{l33} Let $f_0\in \Cf^{Y^\otimes k_0},
 g_j\in \Cf^{H\otimes Y^{\otimes k_j}}$.
Then \be\label{38} E\sum_{l_1,\dots,l_m} (g_0, \bigotimes_{i\in
I_{0}}y_{l_{\sigma(i)}})_{Y^{\otimes
k_0}}\cdot\prod_{j=1}^v\delta\Big( (g_j,\bigotimes_{i\in
I_{j}}y_{l_{\sigma(i)}})_{Y^{\otimes k_j}}\Big)
 \leq \Cs(v,d,\psi)\cdot E\left[ |g_0
|_{v}\prod_{j=1}^v|g_j|_{v-1}\right]. \ee
\end{lem}
\emph{Remark.}  The left hand side of (\ref{372}) can be rewritten
to the form \be\label{380} E\sum_{l_1,\dots,l_m}\delta((g,
y_{l_1})_Y)\delta((g, y_{l_1})_Y)\delta((g, y_{l_2})_Y)\delta((g,
y_{l_2})_Y)\dots \delta((g, y_{l_m})_Y)\delta((g, y_{l_m})_Y). \ee
If $v=2m, k_0=0, k_1=\dots=k_{2m}=1$, $g_0=1,
g_1=\dots=g_{2m}=g\in \Cf^{H\otimes Y}$, then the left hand side
of (\ref{38}) coincides with the expression written in
(\ref{380}). Thus Lemma \ref{l33} implies estimate (\ref{372}).

\emph{Proof of the lemma.} We use induction by $v$. For $v=0,$
conditions imposed  on $\sigma$ can be satisfied   if $m=0$, only
(i.e., if  $g_0$ is a function valued in $\Re$). Thus, for $v=0$,
(\ref{38}) is trivial since $g_0\leq |g_0|\equiv |g_0|_0$. For
$v=1$, conditions imposed on $\sigma$ imply that
$I_0=\{1,\dots,m\}, I_1=\{m+1,\dots,2m\}$ and the function
$\sigma$, restricted to either $I_0$ or $I_1$, is bijective.  Thus
the left hand side of (\ref{38}) can be rewritten to the form
$$
E\sum_{l_1,\dots,l_m} (g_0, \bigotimes_{i=1}^m
y_{l_i})_{Y^{\otimes m}}\cdot\delta\Big( (g_j,\bigotimes_{i=1}^m
y_{l_{\pi(i)}})_{Y^{\otimes m}}\Big),
$$
where $\pi$ is some permutation of $\{1,\dots,m\}$. Using duality
formula (\ref{35}), we rewrite this as
$$
E\Big(Dg_0, A_\pi g\Big)_{H\otimes Y^{\otimes m}},
$$
where the operator $A_\pi\in \Lf(H\otimes Y^{\otimes m})$ is
defined by \be\label{3_trans} A[h\otimes y_{l_1}\otimes
\dots\otimes y_{l_m}]=h\otimes y_{l_{\pi(1)}}\otimes \dots\otimes
y_{l_{\pi(m)}}. \ee
 One can easily see that $A_\pi$ is an isometry
operator, and thus Proposition \ref{p33} provides that (\ref{38})
holds true for $v=1$ with $\Cs(1,d,\psi)=1$.

Suppose that, for some $V\geq 2$, (\ref{38}) holds true for all
$v\leq V-1$. Let us prove that (\ref{38})  holds for $v=V$, also.
For every $l_1,\dots, l_m$, take
$$g=(g_1,\bigotimes_{i\in
I_{1}}y_{l_{\sigma(i)}})_{Y^{\otimes k_1}},\quad f=(g_0,
\bigotimes_{i\in I_{0}}y_{l_{\sigma(i)}})_{Y^{\otimes
k_0}}\cdot\prod_{j=2}^{V}\delta\Big( (g_j,\bigotimes_{i\in
I_{j}}y_{l_{\sigma(i)}})_{Y^{\otimes k_j}}\Big)
$$
and apply duality formula (\ref{35}). Then the left hand side of
(\ref{38}) transforms to the form
$$ E\sum_{l_1,\dots,l_m}\left(D (g_0, \bigotimes_{i\in
I_{0}}y_{l_{\sigma(i)}})_{Y^{\otimes k_0}},(g_1,\bigotimes_{i\in
I_{1}}y_{l_{\sigma(i)}})_{Y^{\otimes k_1}}\right)_H
\cdot\prod_{j=2}^{V}\delta\Big( (g_j,\bigotimes_{i\in
I_{j}}y_{l_{\sigma(i)}})_{Y^{\otimes k_j}}\Big)+ $$
$$
+\sum_{r=2}^{V}E\sum_{l_1,\dots,l_m} (g_0, \bigotimes_{i\in
I_{0}}y_{l_{\sigma(i)}})_{Y^{\otimes
k_0}}\cdot\left(D\left(\delta\Big( (g_r,\bigotimes_{i\in
I_{r}}y_{l_{\sigma(i)}})_{Y^{\otimes
k_r}}\Big)\right),(g_1,\bigotimes_{i\in
I_{1}}y_{l_{\sigma(i)}})_{Y^{\otimes k_1}} \right)_H\times
$$
\be\label{39}
\times\prod_{j\in\{2,\dots,V\}\backslash\{r\}}\delta\Big(
(g_j,\bigotimes_{i\in I_{j}}y_{l_{\sigma(i)}})_{Y^{\otimes
k_j}}\Big). \ee Let us estimate every summand in (\ref{39})
separately. The idea is that every such summand can be written as
\be\label{310} E\sum_{l_1,\dots,l_{\tilde m}} (\tilde g_0,
\bigotimes_{i\in \tilde I_{0}}y_{l_{\tilde
\sigma(i)}})_{Y^{\otimes \tilde k_0}}\cdot\prod_{j=1}^v\delta\Big(
(\tilde g_j,\bigotimes_{i\in \tilde I_{j}}y_{l_{\tilde
\sigma(i)}})_{Y^{\otimes \tilde k_j}}\Big) \ee with $v=V-1$ or
$v=V-2$ and some new $\tilde m, \tilde k_0,\dots, k_{\tilde m} v,
\tilde g_0,\tilde g_{v}, \tilde \sigma$, and thus the inductive
supposition can be applied.

Consider the first summand. Denote by $J$ the set of indices
$r\in\{1,\dots,m\}$ such that $\sigma^{-1}(\{r\})\subset I_0\cup
I_1$. In order to shorten notation, we suppose that
$J=\{1,\dots,\#J\}$ (this does not restrict generality since one
can make an appropriate permutation of the set $\{1,\dots,m\}$ in
order to provide such a property).  Take permutations
$\pi_0:I_0\to I_0$ and $\pi_1:I_1\to I_1$ such that
$$
[\sigma\circ\pi_0](i)=i, \quad i\in\{1,\dots,\#J\}, \qquad
[\sigma\circ\pi_1](i)=i-k_0,\quad i\in\{k_0+1,\dots,k_0+\#J\}.
$$
Then the first summand in (\ref{39}) can be rewritten to the form
$$
E\sum_{l_{\#J+1},\dots,l_m}\left[\sum_{l_1,\dots, l_{\#J}}\left(D
(A_{\pi_0}g_0, \bigotimes_{i=1}^{k_0}
y_{l_{\sigma(\pi_0(i))}})_{Y^{\otimes
k_0}},(A_{\pi_1}g_1,\bigotimes_{i=k_0+1}^{k_0+k_1}
y_{l_{\sigma(\pi_1(i)}})_{Y^{\otimes k_1}}\right)_H\right]\times
$$
$$
\times\prod_{j=2}^{V}\delta\Big( (g_j,\bigotimes_{i\in
I_{j}}y_{l_{\sigma(i)}})_{Y^{\otimes k_j}}\Big),
$$
where the operators $A_{\pi_0}$ and $A_{\pi_1}$ are defined
analogously to (\ref{3_trans}). Denote $$\tilde
g_0=(DA_{\pi_0}g_0, A_{\pi_1} g_1)_{H\otimes Y^{\otimes (\#J)}},$$
then
$$
\left[\sum_{l_1,\dots, l_{\#J}}\left(D (A_{\pi_0}g_0,
\bigotimes_{i=1}^{k_0} y_{l_{\sigma(\pi_0(i))}})_{Y^{\otimes
k_0}},(A_{\pi_1}g_1,\bigotimes_{i=k_0+1}^{k_0+k_1}
y_{l_{\sigma(\pi_1(i)}})_{Y^{\otimes
k_1}}\right)_H\right]=
$$
$$
=\left(\tilde g_0,
\left[\bigotimes_{i=\#J+1}^{k_0}y_{l_{\sigma(\pi_0(i))}}\right]
\otimes\left[\bigotimes_{i=k_0+\#J+1}^{k_1}y_{l_{\sigma(\pi_1(i))}}\right]
\right)_{Y^{k_0+k_1-2\#J}}.
$$
 Put $\tilde
m=m-\#J, \tilde k_0=k_0+k_1-2\#J, \tilde k_1=k_2,\dots, \tilde
k_{V-1}=k_{V}$ and let $\{\tilde I_0,\tilde I_1,\dots,\tilde
I_{V-1}\}$ be the partition of $\{1,\dots, 2\tilde m\}$
corresponding to the family $\{\tilde k_0,\dots, \tilde
k_{V-1}\}$. Put $\tilde g_1=g_2,\dots \tilde g_{V-1}=g_V$ ($\tilde
g_0$ is already defined). At last, define function $\tilde \sigma$
by
$$
\tilde \sigma(i)=\begin{cases}\sigma(\pi_0(i+\#J)),&i=1,\dots,
k_0-\#J\\
\sigma(\pi_1(i+2\#J)),&i= k_0-\#J+1,\dots,k_0+k_1-2\#J\\
\sigma(i+2\#J), &i=k_0+k_1,\dots, 2\tilde m
\end{cases}.
$$
Under such a notation, the first summand in (\ref{39}) has exactly
the form (\ref{310}) with $v=V-1$, and the inductive supposition
provides that this summand is dominated by the term
$$
\Cs(V-1,d,\psi)\cdot E\left[ |\tilde g_0
|_{V-1}\prod_{j=1}^{V-1}|\tilde g_j|_{V-2}\right].
$$
Since $A_{\pi_0},A_{\pi_1}$ are isometric operators, we can apply
Proposition \ref{p33} and obtain that
$$
|\tilde g_0|_{V-1}=|(DA_{\pi_0}g_0, A_{\pi_1} g_1)_{H\otimes
Y^{\otimes (\#J)}}|_{V-1}\leq \Cs(V-1)|DA_{\pi_0}g_0|_{V-1}
|A_{\pi_1} g_1|_{V-1}\leq \Cs(V-1)|g_0|_V\cdot |g_1|_{V-1}.
$$
For every $j=1,\dots, V-1$, $|\tilde
g_j|_{V-2}=|g_{j+1}|_{V-2}\leq |g_{j+1}|_{V-1}.$ Thus, under the
inductive supposition, the first summand in (\ref{39}) is
dominated by the expression given in the right hand side of
(\ref{38}).

All the $V-1$ summands in the second sum in  (\ref{39}) have the
same form and can be estimated similarly; let us make such an
estimation for $r=2$. Using Proposition  \ref{p33}, we rewrite
this summand to the form
$$
E\sum_{l_1,\dots,l_m} (g_0, \bigotimes_{i\in
I_{0}}y_{l_{\sigma(i)}})_{Y^{\otimes k_0}}
\left(B\Big((g_2,\bigotimes_{i\in
I_{r}}y_{l_{\sigma(i)}})_{Y^{\otimes
k_r}}\Big),(g_1,\bigotimes_{i\in
I_{1}}y_{l_{\sigma(i)}})_{Y^{\otimes k_1}} \right)_H
\prod_{j=3}^V\delta\Big( (g_j,\bigotimes_{i\in
I_{j}}y_{l_{\sigma(i)}})_{Y^{\otimes k_j}}\Big)+
$$
$$
+E\sum_{l_1,\dots,l_m} (g_0, \bigotimes_{i\in
I_{0}}y_{l_{\sigma(i)}})_{Y^{\otimes k_0}}\delta\left(\left(D\Big(
(g_2,\bigotimes_{i\in I_{r}}y_{l_{\sigma(i)}})_{Y^{\otimes
k_r}}\Big),(g_1,\bigotimes_{i\in
I_{1}}y_{l_{\sigma(i)}})_{Y^{\otimes k_1}}
\right)_H\right)\times
$$
$$
\times \prod_{j=3}^V\delta\Big( (g_j,\bigotimes_{i\in
I_{j}}y_{l_{\sigma(i)}})_{Y^{\otimes k_j}}\Big)+
$$
$$
 +E\sum_{l_1,\dots,l_m} (g_0, \bigotimes_{i\in
I_{0}}y_{l_{\sigma(i)}})_{Y^{\otimes k_0}} \left(\Big[D\Big(
(g_2,\bigotimes_{i\in I_{r}}y_{l_{\sigma(i)}})_{Y^{\otimes
k_r}}\Big)\Big]^*,(Dg_1,\bigotimes_{i\in
I_{1}}y_{l_{\sigma(i)}})_{Y^{\otimes k_1}} \right)_{H\otimes H}
\times
$$
\be\label{311}
 \times \prod_{j=3}^V\delta\Big( (g_j,\bigotimes_{i\in
I_{j}}y_{l_{\sigma(i)}})_{Y^{\otimes k_j}}\Big). \ee

Let us show that,  after an appropriate rearrangement of the
indices $i$, every summand in (\ref{311}) can be rewritten to the
form (\ref{310}). Such a rearrangement can be organized in the
way, totally analogous to the one used before while the first
summand in (\ref{39}) was estimated. Therefore, in order to
shorten exposition, we do not write here an explicit form for the
permutations of the indices, used in such a rearrangement. The
first summand in (\ref{311}) has the form (\ref{310}) with
$v=V-2$, $\tilde g_j=g_{j+2}, j=1,\dots, V-2$,
$$
\tilde g_0=\Big(A_{\pi_1}[Bg_2], A_{\pi_2} [g_1\otimes
g_0]\Big)_{H\otimes Y^{\otimes\#J^1}},
$$
where $J^1$ is the set of such $i\in \{1,\dots, m\}$  that
$\sigma^{-1}(\{i\})\subset I_0\cup I_1\cup I_2$ (we do not write
here an explicit expressions neither for the permutations
$\pi_1,\pi_2$ nor for $\tilde k_0,\dots, \tilde k_{V-2}, \tilde
\sigma$). Under inductive supposition, this summand is estimated
by
$$
\Cs(V-2,d,\psi)\,E|\tilde g_0|_{V-2}\prod_{j=3}^V |g_j|_{V-3}\leq
\Cs(V,d,\psi)\,E|Bg_2|_{V-2}\cdot|g_1|_{V-2}\cdot|g_0|_{V-2}\prod_{j=3}^V
|g_j|_{V-3}\leq
$$
$$
\leq \Cs(V,d,\psi)
E|g_2|_{V-2}\cdot|g_1|_{V-2}\cdot|g_0|_{V-2}\prod_{j=3}^V
|g_j|_{V-3}\leq \Cs(V,d,\psi) E|g_0|_{V}\prod_{j=1}^V |g_j|_{V-1},
$$
here we used Proposition \ref{p33} and Lemma \ref{l32}.  The
second  summand in (\ref{311}) has the form (\ref{310}) with
$v=V-1$, $\tilde g_j=g_{j+1}, j=2,\dots, V-1$, $\tilde g_0=g_0$,
$$
\tilde g_1=\Big(D[A_{\pi_1}g_2], A_{\pi_2} g_1\Big)_{H\otimes
Y^{\otimes\#J^2}},
$$
where $J^2$ is the set of such $i\in \{1,\dots, m\}$  that
$\sigma^{-1}(\{i\})\subset I_1\cup I_2$. This summand again is
estimated by
$$
\Cs(V-1,d,\psi)\,E|g_0|_{V-1}|\tilde g_1|_{V-2}\prod_{j=3}^V
|g_j|_{V-2}\leq
$$
$$
\leq \Cs(V,d,\psi)\,E |g_0|_{V-1}|
|Dg_2|_{V-2}\cdot|g_1|_{V-2}\cdot|g_0|_{V-2}\prod_{j=3}^V
|g_j|_{V-3}\leq \Cs(V,d,\psi)\, E|g_0|_{V}\prod_{j=1}^V
|g_j|_{V-1}.
$$
At last, the third summand in (\ref{311}) has the form (\ref{310})
with $v=V-2$, $\tilde g_j=g_{j+2}, j=1,\dots, V-2$,
$$
\tilde g_0=\Big(A_{\pi_1}g_0, ([DA_{\pi_2} g_2]^*,
DA\pi_3g_2)_{H\otimes H}\Big)_{Y^{\otimes\#J^1}},
$$
and again is estimated by
$$
\Cs(V-2,d,\psi)\,E|\tilde g_0|_{V-2}|\prod_{j=3}^V |g_j|_{V-3}\leq
$$
$$
\leq \Cs(V,d,\psi)\,E |g_0|_{V-2}|\cdot
|Dg_2|_{V-2}\cdot|Dg_1|_{V-2}\prod_{j=3}^V |g_j|_{V-3}\leq
\Cs(V,d,\psi)\, E|g_0|_{V}\prod_{j=1}^V |g_j|_{V-1}.
$$
The estimates given above show that (\ref{38}) holds for $v=V$ as
soon as it holds  for $v=V-2$ and $v=V-1$. We have already proved
that (\ref{38}) holds for $v=0,1$. Thus, (\ref{38}) holds for
every $v$. The lemma is proved.

\emph{Proof of Lemma \ref{l31}.} We have already proved (\ref{36})
to hold for every $g\in \Cf^H$. Now, let $g\in
W_{2m}^{2m+l-1}(H)$. Consider $\{g_n\}\subset\Cf^H$ such that
$g_n\to g$ in $W_{2m}^{2m+l-1}(H)$ (recall that $\Cf^H$ is dense
in any $W_p^k(H)$ by definition). By (\ref{36}), for any
$k=0,\dots, l$,
$$
\|D^k\delta(g_n)-D^k\delta(g_N)\|_{L_{2m}(\Omega, P, H^{\otimes
k})}\leq \Cs(m,l,d,\psi)\|g_n-g_N\|_{2m,2m+l-1}\to 0, \quad
n,N\to+\infty.
$$
Thus there exist $F_k\in L_{2m}(\Omega, P, H^{\otimes k}),
k=0,\dots, l$ such that
$$
\|D^k\delta(g_n)-F_k\|_{L_{2m}(\Omega, P, H^{\otimes k})}\to
0,\quad n\to +\infty, \quad k=0,\dots, l.
$$
Since operator $\delta$ is closed, $F_0=\delta(g)$. Using that
operator $D$ is closed, one can verify inductively  that
$F_{k}=DF_{k-1}, k=1,\dots, l$. This means that $\delta(g)\in
W_{2m}^l$ with $D^k\delta(g)=F_k, k=0,\dots, l$. At last, using
(\ref{36}) we get
$$
\|\delta(g)\|_{2m,l}^{2m}=E\sum_{k=0}^l \|F_k\|^{2m}_{H^{\otimes
k}}\leq \lim\sup_{n}E\sum_{k=0}^l
\|D^k\delta(g_n)\|^{2m}_{H^{\otimes k}}=
$$
$$
=\lim\sup_{n}\|\delta(g_n)\|_{2m,l}^{2m}\leq \Cs(m,l,d,\psi)
\lim\sup_{n}\|g_n\|_{2m,3m-1}^{2m}=\Cs(m,d,\psi)
\|g\|_{2m,2m+l-1}^{2m}.
$$
The lemma is proved.

\subsection{Malliavin's representation for the densities of the truncated distributions of smooth functionals.}

The typical result in the Malliavin calculus on the Wiener space
is that, when the components $f_1,\dots, f_d$ of a random vector
$f=(f_1,\dots, f_d)$ are smooth enough and the Malliavin matrix
$\sigma^f=\{(Df_i, Df_j)_{H}\}_{i,j=1}^d$ is non-degenerate in a
sense that \be\label{3111}[\det \sigma^f]^{-1}\in \bigcap_{p\geq
1}L_p(\Omega,\Ff,P),\ee  the distribution of $f$ has a smooth
density (see, for instance, \cite{nualart} \S 3.2). Such kind of a
result is useless in the framework, introduced in subsection 3.1,
 since there does not exist any functional $f$ satisfying
(\ref{3111}): if $\eps_1=\dots=\eps_n=0$ then $Df=0$ for every
$f\in\Cf$. In order to overcome this difficulty we use the
following truncation procedure: we consider, instead of $P$, a new
(non-probability) measure $P_\Xi(\cdot)=P(\cdot\cap\Xi)$ with some
set $\Xi\in \sigma (\eps,\zeta)$. If this set is chosen in such a
way that (\ref{3111}) holds true with $P$ replaced by $P_\Xi$ then
the Malliavin's calculus can be applied in order to investigate
the law of $f$ w.r.t. $P_\Xi$. In this subsection, we give the
Malliavin's representation for the density of this law. All
principal steps in our consideration are analogous to those in the
standard Malliavin calculus on the Wiener space (see, for
instance, \cite{nualart}, Chapter 3). Therefore, we sketch  the
proofs only.

Let $f_1,\dots, f_d\in\Cf$ be fixed, consider the Malliavin matrix
$\sigma^f=(\sigma^f_{ij})_{i,j=1}^d$,
$$
\sigma^f_{ij}=(Df_i,
Df_j)_{H}=\sum_{k,r}\psi(\eta_k)[\prt_{kr}f_i((\eta,\eps,\zeta))]\cdot
[\prt_{kr}f_j((\eta,\eps,\zeta))].
$$
Consider a set $\Xi\in\sigma(\eps,\zeta)$ such that $\Xi\subset
\{\det \sigma^f>0\}$ and \be\label{31112}E\1_\Xi[\det
\sigma^f]^{-p}<\infty,\quad p\geq 1.\ee Put $$ \varrho^{f,\Xi}
(\omega)=\begin{cases}[\sigma^f(\omega)]^{-1},& \omega\in \Xi\\
0,&\omega\not\in\Xi
\end{cases}.
$$

\begin{prop}\label{p34} $\varrho^{f,\Xi}\in W_\infty^\infty(\Re^{d\times d})$ and
\be\label{313} (D\varrho^{f,\Xi}, h)_H=-\varrho^{f,\Xi} \cdot
(D\sigma^f, h)_H \cdot \varrho^{f,\Xi},\quad h\in H. \ee
\end{prop}

{\it Sketch of the proof.} It is enough to prove that
$\varrho^{f,\Xi}\in \bigcap_{p\geq 1} W_p^1(\Re^{d\times d})$ and
(\ref{313}) holds true. Suppose that $\sigma^f\geq c I_{\Re^d}$
with some $c>0$. Then one can easily see that $\varrho^{f,\Xi}\in
\Cf$ and (\ref{313}) follows from the well known formula for the
derivative of the inverse matrix,
$$
{d\over dt}[A(t)]^{-1}=-[A(t)]^{-1}\cdot[{d\over dt}A(t)]\cdot
[A(t)]^{-1}.
$$
%valid at the point $t=t_0$ as soon as the function $A\Re\ni
%t\mapsto A(t)\in\Re^{d\times d}$ is differentiable at the point
%$t_0$ and $A(t_0)$ is invertible.
In the general case, consider the matrix-valued functions
$\sigma^{f,c}=\sigma^f+cI_{\Re^d}$ and
$\varrho^{f,\Xi,c}=\1_\Xi\cdot [\sigma^{f,c}]^{-1}$, $c>0$.
Condition (\ref{31112}) provides that $\varrho^{f,\Xi,c}\to
\varrho^{f,\Xi}, c\to 0+$ in any $L_p$. It is already proved that
(\ref{313}) holds true for the functionals indexed by $c$. Thus,
passing to the limit as $c\to 0+$, we obtain the required
statement.

Denote $\vartheta_i^{f,\Xi}=\sum_{k=1}^d \varrho_{ki}^{f,\Xi}\cdot
Df_k,i=1,\dots, d$. Also denote, by $E_\Xi$, the expectation
w.r.t. $P_\Xi$.

\begin{prop}\label{p35} For every  $i=1,\dots, d$, $
\upsilon\in W_\infty^\infty$ and every $F\in C_b^\infty(\Re^d)$,
\be\label{314} E_\Xi[\prt_i F](f_1,\dots, f_d)\cdot
 \upsilon=E_\Xi F(f_1,\dots, f_d)\cdot
\delta\left(\upsilon\cdot\vartheta_i^{f,\Xi} \right). \ee
\end{prop}

{\it Sketch of the proof.}  It follows from Propositions
\ref{p33}, \ref{p34} and Lemma \ref{l31} that $\upsilon \cdot
\vartheta_i^f \in Dom(\delta)$. Since $\Xi\in \sigma
(\eps,\zeta)$, the function $\1_\Xi$ belongs to $\Cf$ and has its
stochastic derivative equal to 0. Proposition \ref{p32} provides
that $\delta(\1_{\Xi}\cdot g)=\1_\Xi \delta(g), g\in Dom(\delta)$.
Therefore
 $$
 E[\prt_i F](f_1,\dots, f_d)\cdot \1_{\Xi}\cdot \upsilon
 =\sum_{j=1}^d E[\prt_j F](f_1,\dots, f_d)\cdot
\1_{\Xi}\cdot
 \upsilon\cdot [\sigma^f \cdot \varrho^{f,\Xi}]_{ij}=
 $$
 $$
  =\sum_{j=1}^d \sum_{k=1}^d E[\prt_j F](f_1,\dots, f_d)\cdot
\1_{\Xi}\cdot
 \upsilon\cdot \varrho_{ki}^{f,\Xi}\cdot\sigma^f_{jk}=
$$
$$
=E\sum_{k=1}^d \1_{\Xi}\cdot
 \upsilon\cdot \varrho_{ki}^{f,\Xi}\left(\sum_{j=1}^d [\prt_j F](f_1,\dots, f_d)\cdot D
 f_j, Df_k\right)_H=E\left(D[F(f_1,\dots, f_d)],\1_{\Xi}\cdot
 \upsilon\cdot\vartheta_i^{f,\Xi} \right)_H=
$$
$$
=EF(f_1,\dots, f_d)\cdot\delta\left(\1_{\Xi}\cdot
 \upsilon \cdot \vartheta_i^f \right)=EF(f_1,\dots, f_d)\cdot\1_{\Xi}\cdot \delta\left(
 \upsilon \cdot \vartheta_i^f \right),
$$
that provides (\ref{314}).

Put
$$
\upsilon_1^{f,\Xi}=1,\quad \upsilon_{l+1}^{f,\Xi}=\delta\left(
 \upsilon_l^{f,\Xi} \cdot \vartheta_l^{f,\Xi} \right), \quad l=1,\dots, d,
$$
$$
\Upsilon^{f,\Xi}=\upsilon_{d+1}^{f,\Xi},\quad
\Upsilon^{f,\Xi}_{i}=\delta\left(
 \Upsilon^{f,\Xi} \cdot \vartheta_i^{f,\Xi} \right),\quad i=1,\dots,d.
 $$
Write $P^f_\Xi$ for the distribution of $f$ w.r.t. $P_\Xi$.  For
$\alpha_1,\dots,\alpha_d\in\{0,1\}$, denote
$$
\1_{\alpha_1\dots\alpha_d}(x)=\1_{(-1)^{\alpha_1}x_1\geq 0,\dots,
(-1)^{\alpha_d}x_d\geq 0},\quad x\in\Re^d.
$$

\begin{prop}\label{p36} The distribution $P_\Xi^f$ has a
 density $p^f_\Xi$, bounded together with all its derivatives
 $\prt_ip_\Xi^f, i=1,\dots,d$. For any
$\alpha_1,\dots,\alpha_d\in\{0,1\}$, \be\label{315}
p_\Xi^f(y)=(-1)^{\alpha_1+\dots+
\alpha_d}E\1_{\alpha_1\dots\alpha_d}(f-y)
\cdot\Upsilon^{f,\Xi},\quad y\in \Re^d,\ee \be\label{316}
\prt_ip_\Xi^f(y)=(-1)^{\alpha_1+\dots+
\alpha_d+\alpha_i+1}E\1_{\alpha_1\dots\alpha_d}(f-y)
\cdot\Upsilon^{f,\Xi}_{i},\quad y\in \Re^d.\ee
\end{prop}

{\it Sketch of the proof.} Applying iteratively (\ref{314}) one
can deduce that, for every $F\in C_b^\infty(\Re^d)$,
\be\label{3141} E_\Xi[\prt_1\dots\prt_d F](f_1,\dots, f_d)=E
F(f_1,\dots, f_d)\cdot \Upsilon^{f,\Xi},\ee \be\label{3142}
E_\Xi[\prt_1\dots\prt_d\prt_i F](f_1,\dots,f_d)=E F(f_1,\dots,
f_d)\cdot \Upsilon^{f,\Xi}_{i},\quad i=1,\dots,d. \ee Now, the
informal way to get representation (\ref{315}) is to apply
(\ref{3141}) to $F=\1_{\alpha_1\dots\alpha_d}$:
$$
p_\Xi^f(y)=(-1)^{\alpha_1+\dots+\alpha_d+d}\prt_1\dots\prt_d
E_\Xi\1_{\alpha_1\dots\alpha_d}(f-y)=(-1)^{\alpha_1+\dots+\alpha_d}
E_\Xi[\prt_1\dots\prt_d\1_{\alpha_1\dots\alpha_d}](f-y)=
$$
\be\label{317} =(-1)^{\alpha_1+\dots+\alpha_d}
E\1_{\alpha_1\dots\alpha_d}(f-y) \cdot\Upsilon^{f,\Xi} \ee In
order to justify (\ref{317}) one should consider smooth
approximations $F_n$ for the function
$F=\1_{\alpha_1\dots\alpha_d}$ and use Fubini theorem (we omit
detailed exposition here, referring the reader, for instance, to
\cite{nualart}, \S\S 3.1, 3.2). Similarly, (\ref{316}) is provided
by (\ref{3142}) and the formula
$$
\prt_ip_\Xi^f(y)=(-1)^{\alpha_1+\dots+\alpha_d+1}
E_\Xi[\prt_1\dots\prt_d\prt_i\1_{\alpha_1\dots\alpha_d}](f-y).
$$

\subsection{Estimates for the densities of the truncated distributions of smooth functionals.}

Proposition \ref{p36} immediately provides the following family of
estimates for the density $p_\Xi^f$ of the truncated distribution
of $f$.

\begin{cor}\label{c31} For any $y\in\Re^d$,
$$
p^f_\Xi(y)\leq
\|\Upsilon^{f,\Xi}\|_{L_2}\cdot\min_{\alpha_1,\dots,\alpha_d\in\{0,1\}}P^{1\over
2}_\Xi\Big((-1)^{\alpha_1}f_1\geq (-1)^{\alpha_1}y_1,\dots,
(-1)^{\alpha_d}f_d\geq (-1)^{\alpha_d} y_d\Big)\leq
\|\Upsilon^{f,\Xi}\|_{L_2},
$$
$$
\prt_ip^f_\Xi(y)\leq
\|\Upsilon_{i}^{f,\Xi}\|_{L_2}\cdot\min_{\alpha_1,\dots,\alpha_d\in\{0,1\}}P^{1\over
2}_\Xi\Big((-1)^{\alpha_1}f_1\geq (-1)^{\alpha_1}y_1,\dots,
(-1)^{\alpha_d}f_d\geq (-1)^{\alpha_d} y_d\Big)\leq
\|\Upsilon_{i}^{f,\Xi}\|_{L_2},
$$
$i=1,\dots, d$. In particular, $p_\Xi^f$ satisfies Lipschitz
condition with the constant
$L=\sum_{i=1}^d\|\Upsilon_{i}^{f,\Xi}\|_{L_2}$.
\end{cor}

In this subsection we give explicit estimates for
$\|\Upsilon^{f,\Xi}\|_{L_2},\|\Upsilon_{i}^{f,\Xi}\|_{L_2},
i=1,\dots,d$. Our estimates somewhat differ from the standard
Malliavin-type ones. In our considerations, we operate with the
matrix $[\sigma^f]^{-1}$ straightforwardly and do not use (unlike
in the standard Maliavin's approach) representation of this matrix
via the Cramer's formula $[\sigma^f]^{-1}=[\det
\sigma^f]^{-1}\cdot \Sigma^f$ ($\Sigma^{f}$ denotes the cofactor
matrix for $\sigma^f$). This is caused by our goal to prove,
together with existence of the density, an explicit estimates for
it like the estimate (ii) of  Theorem \ref{t11}.

Let us give an iterative description of the family
$\{\upsilon_l^{f,\Xi}\}$ involved into construction of
$\Upsilon^{f,\Xi},\Upsilon_{i}^{f,\Xi}$.  We introduce two
families of operators acting on $W_\infty^\infty$:
$$
I_{i}:\phi\mapsto \delta(\phi Df_i), \quad J_{ijk}:\phi\mapsto
\phi(D\sigma^f_{jk}, Df_i)_H, \quad i,j,k=1,\dots, d.
$$
 We call any operator $I_1,\dots,
 I_d$ an \emph{operator of the type $I$},  and any
 operator from the set $\{J_{ijk},i,j,k=1,\dots, d\}$ an \emph{operator of the type $J$}.
 We denote by $\Kf(m,M)$ the class of all functions that can be
obtained from $\phi\equiv 1$ by applying,   in arbitrary order, of
$m$ operators of the type $I$ and $M$ operators of the type $J$.

\begin{prop}\label{p37} For any $l=1,\dots, d+1$, there exist constant $\Cs(d,l)\in\NN$ such
that $\upsilon_l^{f,\Xi}$ is a sum of at most $\Cs(d,l)$ summands
of the type \be\label{318}
\phi\cdot\prod_{k=1}^r\varrho^{f,\Xi}_{i_k j_k}, \ee where
$i_k,j_k=1,\dots, d$ are arbitrary and $\phi$ belongs to some
class $\Kf(m,M)$ with $m+M=l-1$ and $r=M+l-1$.
\end{prop}

\demo  We use induction by $l$. For $l=1$, the statement is
trivial since $\upsilon_1^{f,\Xi}=1\in \Kf(0,0).$ Suppose the
statement of the Lemma to hold true for some $l\leq d$. Let us
prove this statement  for $l+1$. Due to the inductive supposition,
$\upsilon_{l}^{f,\Xi}$ is a sum of at most $\Cs(d,l)$ summands of
the type $$ \delta \Big(\phi\cdot\prod_{k=1}^r\varrho^{f,\Xi}_{i_k
j_k}\cdot\vartheta_{l}^{f,\Xi}\Big), $$ $\phi\in\Kf(m,M), m+M=l,
M+l=r$. We have
$$
\delta \Big(\phi\cdot\prod_{k=1}^r\varrho^{f,\Xi}_{i_k
j_k}\vartheta_{l}^{f,\Xi}\Big)=\sum_{q=1}^d \delta
\Big(\phi\cdot[\prod_{k=1}^r \varrho^{f,\Xi}_{i_k j_k}]\cdot
\varrho_{ql}^{f,\Xi}\cdot Df_l\Big)
$$
Thus, $\upsilon_{l}^{f,\Xi}$ is a sum of at most $d\cdot \Cs(d,l)$
summands of the type $$ \delta
\Big(\phi\cdot\prod_{k=1}^{r+1}\varrho^{f,\Xi}_{i_k j_k}\cdot
Df_l\Big), $$ $\phi\in\Kf(m,M), m+M=l, M+l=r$. Due to Propositions
\ref{p32} and \ref{p34},
$$ \delta
\Big(\phi\cdot\prod_{k=1}^{r+1}\varrho^{f,\Xi}_{i_k j_k}\cdot
Df_l\Big)= \delta(\phi\cdot Df_l\Big)\cdot
\prod_{k=1}^{r+1}\varrho^{f,\Xi}_{i_k
j_k}-\left(D\prod_{k=1}^{r+1}\varrho^{f,\Xi}_{i_k j_k}, \phi\cdot
Df_l\right)_H=
$$
\be\label{319}=\delta(\phi\cdot Df_l)\cdot
\prod_{k=1}^{r+1}\varrho^{f,\Xi}_{i_k
j_k}-\sum_{q=1}^{l+1}\phi\cdot\left[\prod_{k\leq r+1,
k\not=q}\varrho^{f,\Xi}_{i_k
j_k}\right]\Big(\varrho^{f,\Xi}\cdot(D\sigma^f,
Df_l)\cdot\varrho^{f,\Xi}\Big)_{i_qj_q}. \ee Since $\phi\in
\Kf(m,M)$, $\delta(\phi \cdot Df_l)\in \Kf (m+1, M)$. Thus, the
first term in the right hand side of (\ref{319}) has the form
(\ref{318}). Every  summand in the sum in the right hand side of
(\ref{319}) is a sum of $d^2$ terms of the type
$$
\phi\cdot (D\sigma^f_{\tilde i\tilde j},
Df_l)_H\cdot\left[\prod_{k=1}^{r+2}\varrho^{f,\Xi}_{\tilde i_k
\tilde j_k}\right]
$$
Every such a term has the form (\ref{318}), since $\phi\cdot
(D\sigma^f_{\tilde i\tilde j}, Df_l)_H\in \Kf(m, M+1)$  for
$\phi\in\Kf(m, M)$. Therefore, the statement of the Lemma holds
true for $l+1$, also, with $\Cs(d,l+1)=\Cs(d,l)[1+d^2(l+1)]$. The
proposition is proved.

Recall that $\Upsilon^{f,\Xi}=\upsilon^{f,\Xi_{d+1}}$, and thus
Proposition \ref{p37} provides that $\Upsilon^{f,\Xi}$ is a sum of
not more than $\Cs(d)$ summands of the type (\ref{318}) with
$\phi\in \Kf(d-M, M)$ and $r=M+d$ ($M$ may vary from $0$ to $d$).
For every such a summand,
\be\label{332}\|\phi\cdot\prod_{k=1}^{M+d}\varrho^{f,\Xi}_{i_k
j_k}\|_{L_2}\leq \|\phi\|_{L_4}\cdot
\Big[E\|\varrho^{f,\Xi}\|_{\mathcal{M}}^{4(M+d)}\Big]^{{1\over
4}}, \ee where $\|A\|_{\mathcal{M}}=\max_{i,j=1,\dots, d}|A_{ij}|,
A\in \Re^{d\times d}$. Thus, in order to estimate
$\|\Upsilon^{f,\Xi}\|_{L_2}$, it is sufficient to estimate
$\max_{\phi\in \Kf(d-M, M)}\|\phi\|_{L_4}$. Denote $\alpha_i=Df_i,
\beta_{ijk}=(D\sigma^f_{jk}, Df_i)_H$.

\begin{prop}\label{p38} For every $\phi\in\Kf(d-M,M),$
\be\label{333} \|\phi\|_{L_4}\leq \Cs(d,\psi)\cdot
\Big(\max_i\|\alpha_i\|_{2(d+1)(d+2), (d+1)^2-1}\Big)^{d-M}\cdot
\Big(\max_{ijk}\|\beta_{ijk}\|_{2(d+1)(d+2), (d+1)^2-1}\Big)^{M}.
\ee
\end{prop}
\demo By Lemma \ref{l31} and Proposition \ref{p33}, for any
$i=1,\dots, d, \phi\in W_\infty^\infty, m\geq 0$ $$ \|\delta(\phi
\alpha_i)\|_{2m+2, m^2-1}\leq \Cs(m,d,\psi)\|\phi
\alpha_i\|_{2m+2, (m+1)^2-1}\leq
$$
\be\label{334} \leq  \Cs(m,d,\psi)\|\phi\|_{2m+4,
(m+1)^2-1}\|\alpha_i\|_{2(m+1)(m+2), (m+1)^2-1}.
 \ee
 By Proposition \ref{p33}, for any
$i,j,k=1,\dots, d, \phi\in W_\infty^\infty, m\geq 0$
 $$
 \|\phi \beta_{ijk}\|_{2m+2, m^2-1}\leq \Cs(m)
\|\phi\|_{2m+4, m^2-1} \beta_{ijk}\|_{2(m+1)(m+2), m^2-1}\leq
$$
\be\label{335}\leq \Cs(m,d,\psi)\|\phi\|_{2m+4,
(m+1)^2-1}\|\beta_{ijk}\|_{2(m+1)(m+2), (m+1)^2-1}.
 \ee
Recall that $\|\cdot \|_{L_4}=\|\cdot \|_{4, 0}$ and
$\phi\in\Kf(d-M, M)$ is obtained from $1$ by applying (in some
order) of $d-M$ operators of the type $I$ and $M$ operators of the
type $J$. Thus, in order to obtain (\ref{333}), one
 should consequently put $m=d,d-1,\dots, 1$ and apply
  either inequality (\ref{334}) or inequality (\ref{335})
 depending on what type of the operator ($I$ or $J$) was
 applied at this position  in the  construction of the function $\phi$.
The proposition is proved.

Recall that $\sigma^f_{jk}=(Df_j, Df_k)_H$. By Proposition
\ref{p33},
$$
|\beta_{ijk}|_m\leq \Cs(m)|Df_i|_m|D\sigma^f_{jk} |_m\leq \tilde
\Cs(m)|Df_i|_m |Df_j|_{m+1}|Df_k|_{m+1}.
$$
Thus, Proposition \ref{p38} provides the following estimate for
$\|\Upsilon^{f,\Xi}\|_{L_2}$. Denote
$$
N_d(f)=\max_{i=1,\dots,d}\max_{m=1, \dots, (d+1)^2}\Big[E\|D^m
f_i\|_{H^{\otimes m}}^{2(d+1)(d+2)}\Big]^{1\over 2(d+1)(d+2)}.
$$

\begin{cor}\label{c32} There exists a constant $\Ls_d$, dependent on $d$ and $\psi$ only, such
that \be\label{336} \|\Upsilon^{f,\Xi}\|_{L_2}\leq
\Ls_d\sum_{M=0}^d \Big[N_d(f)\Big]^{d+2M}\cdot
\Big[E\|\varrho^{f,\Xi}\|_{\mathcal{M}}^{4(M+d)}\Big]^{{1\over
4}}.  \ee
\end{cor}

For $\|\Upsilon^{f,\Xi}_i\|_{L_2}$,  the following estimate holds
true (we omit the proof since it is totally analogous to the proof
of (\ref{336}) given above).

\begin{prop}\label{p39} For any $i=1,\dots, d$,
\be\label{337} \|\Upsilon^{f,\Xi}_i\|_{L_2}\leq
\Ls_{d+1}\sum_{M=0}^{d+1} \Big[N_{d+1}(f)\Big]^{d+2M+1}\cdot
\Big[E\|\varrho^{f,\Xi}\|_{\mathcal{M}}^{4(M+d+1)}\Big]^{{1\over
4}}.  \ee \end{prop}

At the end of this section, we formulate a general local limit
theorem. This theorem is a straightforward corollary of the
representation given by  Proposition \ref{p36} and the estimates
(\ref{336}) and (\ref{337}).  Consider the sequence of probability
spaces $\{(\Omega^n,\Ff^n,P^n), n\geq 1\}$ of the type
(\ref{3200}) with the given measures $\pi_U,\kap,\nu$. Let the
function $\psi$ and the functions $\theta_n$ be fixed. Denote
$H_n=\Re^{d\times n}$ and consider the derivative, gradient and
Sobolev spaces constructed in subsection 3.1.  For a sequence of
random vectors $f^n:\Omega^n\to \Re^d$ and a sequence of sets
$\{\Xi_n\in \sigma (\eps, \zeta)\}$ denote by $P^{n,f^n}$ the
distribution of $f^n$ w.r.t. $P^n$ and by $P^{n,f^n}_{\Xi_n}$ the
distribution of $f^n$ w.r.t. $P^n_{\Xi_n}\equiv
P^n(\cdot\cap\Xi_n)$. Denote
$$
\Ks_d(f,\Xi)=\Ls_d\cdot\sum_{M=0}^d \Big[N_d(f)\Big]^{d+2M}\cdot
\Big[E\|\varrho^{f,\Xi}\|_{\mathcal{M}}^{4(M+d)}\Big]^{{1\over
4}}.
$$

\begin{thm}\label{t31} Suppose that $\{f^n\}$ and $\{\Xi_n\}$ satisfy condition

\begin{itemize}
\item[(C1)] $f^n\in W_{2(d+2)(d+3)}^{(d+2)^2}(\Re^d)$ and $\sup_n
\Ks_{d+1}(f^n,\Xi_n)<+\infty$.

\end{itemize}

Then $P^{n,f^n}_{\Xi_n}$ possess a densities $p^{f^n}_{\Xi_n}$.
Moreover,
\begin{itemize}

\item[(a)] $p^{f^n}_{\Xi_n}(y)\leq \Ks_{d}(f^n,\Xi_n)\cdot
P_{\Xi_n}^{1\over 2}(\|f^n\|\geq \|y\|)$;

\item[(b)] $p^{f^n}_{\Xi_n}$ satisfy Lipschitz condition with the
common constant equal to $d\cdot\sup_n \Ks_{d+1}(f^n,\Xi_n)$.

\end{itemize}

\noindent If, additionally,
\begin{itemize}

\item[(C2)] $f^n$ converge in distribution to some random vector
$f$;

\item[(C3)] $P^n(\Xi_n)\to 1, n\to +\infty$;
\end{itemize}

\noindent then  the distribution of the vector $f$ possess a
density $p^f$ and

\begin{itemize}
\item[(c)] $\sup_{y\in\Re^d} |p^{f^n}_{\Xi_n}(y)-p^f(y)|\to 0,
\quad n\to +\infty$.
\end{itemize}

\end{thm}

\demo  Although the statements of Corollaries \ref{c31} and
\ref{c32} are formulated for $f\in \Cf^{\Re^d}$, they can be
extended to $f\in W_{2(d+2)(d+3)}^{(d+2)^2}(\Re^d)$ by a standard
approximation procedure. For any $y\in \Re^d$, there exist a
choice of the signs $\alpha_1,\dots, \alpha_d$ such that
$$
P_{\Xi_n}\Big((-1)^{\alpha_1}f_1^n\geq (-1)^{\alpha_1}y_1,\dots,
(-1)^{\alpha_d}f_d^n\geq (-1)^{\alpha_d} y_d\Big)\leq
P_{\Xi_n}(\|f^n\|\geq \|y\|).
$$
Therefore, statements (a) and (b) follow immediately from
Corollaries \ref{c31} and \ref{c32}. Statement (b) provides that
the sequence $\{p_{\Xi_n}^{f^n}\}$ has a compact closure in the
space $C(\Re^d)$ with the topology of uniform convergence on a
compacts. This together with the conditions (C2),(C3) provides
(c). The theorem is proved.

\section{Proofs of Theorems \ref{t11} -- \ref{t21}}

\subsection{Proof of Theorem \ref{t11}}

We reduce the proof of Theorem \ref{t11} to the verification of
the conditions of Theorem \ref{t31} and explicit estimation of the
expression in the right hand side of (a). We use, without
additional discussion, notation introduced in Section 3.

 Denote $f^n_{x,t}=X_n(t)-x$, where
  the processes $X_n$ are defined by (\ref{01}),(\ref{02}) with the initial value $X_n(0)=x\in\Re^d$.
When it does not cause misunderstanding, we omit the indices $x,t$
and write $f^n$ for $f^n_{x,t}$.

We conduct the proof in several steps. First, we give explicit
expressions for the derivatives of the functionals $f^n$. Next, we
estimate the moments of these derivatives (this allows us to
estimate $N_{d+1}(f^n)$). Then, on the properly chosen $\Xi_n$, we
estimate the inverse matrix for the Malliavin matrix
$\sigma^{f^n}$ (this allows us to estimate
$\Ks_{d+1}(f^n,\Xi_n)$). At last, we estimate the tail
probabilities $P_{\Xi_n}(\|f^n\|\geq \|y\|)$ in order to provide
the estimates given in the statement (ii) of Theorem \ref{t11}.

Everywhere below we suppose conditions (B1), (B2$_\kappa$), (B3)
of Theorem \ref{t11} to hold true. We prove (i) -- (iii) in
details and give a brief sketch of changes that should be made in
order to prove (ii$'$),(iii$'$). We put \be\label{400}
 \theta_n(\zeta)=\1_{\max_{k\leq n}\|\zeta_k\|\leq n^\varsigma},\quad
 \varsigma={\kappa-1\over 2\kappa+2}.
 \ee
 In
 order to make notation more convenient, we rewrite (\ref{01}) to
 the form
\be\label{40}X_n\brakes{k\over n}=X_n\brakes{k-1\over
n}+a\brakes{X_n\brakes{k-1\over n}}\cdot {1\over n} +\sum_{r=1}^d
b_r\brakes{X_n\brakes{k-1\over n}}\cdot {\xi_{kr}\over\sqrt n},
\ee here $\xi_{k1},\dots,\xi_{kd}$ are the components of the
vector $\xi_k$ and $b_1,\dots,b_d$ are the columns of the matrix
$b$.

 \begin{lem}\label{l41} For every $t\in[0,1]$, $X_n(t)\in \bigcap_{p>1}W^{(d+2)^2}_p(\Re^d)$.
 Derivatives  $Y_n(t)=DX_n(t), t\in[0,1]$ satisfy relations
$$ Y_n(0)=0,\quad Y_n\brakes{k\over
n}=Y_n\brakes{k-1\over n}+\nabla a\brakes{X_n\brakes{k-1\over
n}}Y_n\brakes{k-1\over n}\cdot {1\over n} +
$$
\be\label{41} +\sum_{r=1}^d \left[\nabla
b_r\brakes{X_n\brakes{k-1\over n}}Y_n\brakes{k-1\over n}\cdot
{\xi_{kr}\over\sqrt
n}+{\theta_n(\zeta)\1_{\eps_k=1}\psi(\eta_k)\over\sqrt
n}b_r\brakes{X_n\brakes{k-1\over n}}\otimes e_{kr}\right], \ee
 \be\label{42} Y_n(t)=Y_n\brakes{k-1\over
n}+(nt-k+1)\left[Y_n\brakes{k\over n}-Y_n\brakes{k-1\over n}
\right],\quad t\in\left[{k-1\over n},{k\over n}\right),\ee
$k=1,\dots,n.$

\end{lem}

{\it Sketch of the proof.} The proof is quite standard, and thus
we just outline its main steps. Using induction by $k$, one can
easily verify that, for every $j,k,r$, there exists
$\prt_{jr}X_n\brakes{k\over n}=(Y_n\brakes{k\over n},
e_{jr})_{H_n}$ with $Y_n$ defined by (\ref{41}). One can see that
$Y_n\equiv 0$ as soon as $\max_{k\leq n}\|\zeta_k\|>
n^{\varsigma}$ and, therefore,
 $\mathop{\mathrm{ess}\sup}\|Y_n\brakes{k\over
n}\|<+\infty$ for every $k\leq n$. Iterating these considerations,
one can verify that  $\mathop{\mathrm{ess}\sup} \|\nabla_\eta^m
X_n\brakes{k\over n}\|<+\infty$ for every $k\leq n,m\leq (d+2)^2$,
that means that $X_n\brakes{k\over
n}\in\bigcap_{p>1}W^{(d+2)^2}_p(\Re^d)$ with $DX_n\brakes{k\over
n}=Y_n\brakes{k\over n}$, that gives the statement of the Lemma
for $t={k\over n}$. For arbitrary $t\in[0,1]$, this statement
holds by linearity.

Denote $\mu_\kappa(\xi)=E\|\xi_1\|^\kappa$.

 \begin{lem}\label{l42} For every $p\geq 1, m\in \NN$, there exist constant
 $\Cs(a,b,d,U,\mu_\kappa(\xi),m,p)$ such that
 \be\label{43}E\|D^mX_n(t)\|^p_{H^{\otimes m}\otimes \Re^d}\leq
 \Cs(a,b,d,U,\mu_\kappa(\xi),m,p)\cdot t^{{p\over 2}}, \quad t\in \left[{1\over n}, 1\right].
 \ee
  \end{lem}
\demo Consider first the case $m=1$. It is enough to prove
(\ref{43}) for $t={k\over n}, k=1,\dots, n$ and $p=2q, q\in \NN$.
We have
$$
\left\|Y_n\brakes{k\over n}\right\|^2_{H\otimes
\Re^d}=\left\|Y_n\brakes{k-1\over n}\right\|^2_{H\otimes
\Re^d}+{1\over n}\left(Y_n\brakes{k-1\over n},\nabla
a\brakes{X_n\brakes{k-1\over n}}Y_n\brakes{k-1\over n}
\right)_{H\otimes \Re^d}+
$$
$$
+\sum_{r=1}^d\left(Y_n\brakes{k-1\over n},\nabla
b_r\brakes{X_n\brakes{k-1\over n}}Y_n\brakes{k-1\over n}
\right)_{H\otimes \Re^d}\cdot {\xi_{kr}\over\sqrt n}+
$$
$$
+{1\over n^2}\left\|\nabla a\brakes{X_n\brakes{k-1\over
n}}Y_n\brakes{k-1\over n} \right\|_{H\otimes
\Re^d}^2+
$$
$$
+\sum_{r=1}^d\left(\nabla a\brakes{X_n\brakes{k-1\over
n}}Y_n\brakes{k-1\over n},\nabla b_r\brakes{X_n\brakes{k-1\over
n}}Y_n\brakes{k-1\over n} \right)_{H\otimes \Re^d}\cdot
{\xi_{kr}\over n^{3\over 2}}+
$$
$$
+\sum_{r_1,r_2=1}^d\left(\nabla b_{r_1}\brakes{X_n\brakes{k-1\over
n}}Y_n\brakes{k-1\over n},\nabla
b_{r_2}\brakes{X_n\brakes{k-1\over n}}Y_n\brakes{k-1\over n}
\right)_{H\otimes \Re^d}\cdot {\xi_{kr_1}\xi_{kr_2}\over n}+
$$
\be\label{44} +{\theta_n(\zeta)\1_{\eps_k=1}\psi(\eta_k)\over
n}\sum_{r=1}^d\left\|b_r\brakes{X_n\brakes{k-1\over
n}}\right\|^2_{\Re^d}, \ee here we have used the fact that
$Y_n\brakes{k-1\over n}, \nabla a\brakes{X_n\brakes{k-1\over
n}}Y_n\brakes{k-1\over n},\nabla b_r\brakes{X_n\brakes{k-1\over
n}}Y_n\brakes{k-1\over n}, r=1,\dots, d$ belong to the subspace
generated by the vectors of the type $v\otimes e_{jr}, v\in\Re^d,
j<k, r=1,\dots, d$ and $b_r\brakes{X_n\brakes{k-1\over n}}\otimes
e_{kr} ,r=1,\dots, d$ are orthogonal to this subspace.

 Denote
$\left\|Y_n\brakes{k\over n}\right\|^2_{H\otimes
\Re^d}=\Upsilon_k$. Recall that $Y_n\brakes{k\over n}=0,
k=1,\dots, n$ as soon as there exist $j=1,\dots, d$ such that
$\|\zeta_j\|>  n^{\varsigma}$. Since coefficients $a, b$ are
 bounded together with their derivatives, we can rewrite
 (\ref{44}) as
$$
\Upsilon_k=\left[\Upsilon_{k-1}+\Theta_{k-1}\cdot {1\over n} +
\sum_{r=1}^d\Lambda_{k-1,r}\cdot {\xi_{kr}\over\sqrt n}\cdot
\1_{\|\zeta_{k}\|\leq
n^{\varsigma}}+\sum_{r_1,r_2=1}^d\Delta_{k-1,r_1,r_2}\cdot
{\xi_{kr_1}\xi_{kr_2}\over n}\cdot \1_{\|\zeta_{k}\|\leq
n^{\varsigma}}\right]\times
$$
$$
\times \1_{\max_{j<k}\|\zeta_{j}\|\leq n^{\varsigma}},\quad
k=1,\dots, n
$$
with $\Ff_{k-1}\equiv\sigma(\eta_1,\eps_1,\zeta_1,\dots,
\eta_{k-1},\eps_{k-1},\zeta_{k-1})$ -- measurable $\Theta_{k-1},
\Lambda_{k-1,r_1,r_2}$ such that \be\label{45} |\Theta_{k-1}|\leq
\Cs(a,b,d)\left[1+\Upsilon_{k-1}\right], \quad
|\Lambda_{k-1,r_1,r_2}|\leq \Cs(a,b,d)\Upsilon_{k-1}. \ee Since
$\Upsilon_k\geq 0$, we have \be\label{46} E\Upsilon_k^q\leq
E\left(\Upsilon_{k-1}+\Theta_{k-1}\cdot {1\over n}\right)^q+
\sum_{i=0}^{q-1} {q!\over
i!(q-i)!}E\left(\Upsilon_{k-1}+\Theta_{k-1}\cdot {1\over
n}\right)^i\times
$$
$$
 \times\left[\sum_{r=1}^d\Lambda_{k-1,r}\cdot
{\xi_{kr}\over\sqrt n}\cdot \1_{\|\zeta_{k}\|\leq
n^{\varsigma}}+\sum_{r_1,r_2=1}^d\Delta_{k-1,r_1,r_2}\cdot
{\xi_{kr_1}\xi_{kr_2}\over n}\cdot \1_{\|\zeta_{k}\|\leq
n^{\varsigma}}\right]^{q-i}. \ee We have
$\xi_{kr}=\eta_{kr}+\zeta_{kr}$ and the set $U$ of the possible
values of $\eta_{k}=(\eta_{k1},\dots,\eta_{kd})$ is bounded. In
addition, $|\zeta_{kr}|\1_{\|\zeta_{k}\|\leq  n^{\varsigma}}\leq
n^\varsigma$. Therefore, \be\label{450} E\left|{\xi_{kr}\over
\sqrt {n}}\right|^l\cdot \1_{\|\zeta_{k}\|\leq  n^{\varsigma}}
\leq {(\Cs(U)\cdot n^\varsigma)^l\over n^{l\over 2}}= \Cs^l(U)
n^{l(\varsigma-{1\over 2})}\leq {\Cs^l(U)\over n}, \quad l\geq
\kappa+1. \ee Since $E\|\xi_k\|^\kappa<+\infty$, \be\label{451}
E\left|{\xi_{kr}\over \sqrt {n}}\right|^l\cdot
\1_{\|\zeta_{k}\|\leq  n^{\varsigma}} \leq {E|\xi_{kr}|^l\over
n^{l\over 2}}\leq {E(\|\xi_k\|^\kappa\vee 1)\over n}, \quad
l=2,\dots,\kappa. \ee Similarly,
\be\label{4510}\left|E{\xi_{kr_1}\xi_{kr_2}\over {n}}\cdot
\1_{\|\zeta_{k}\|\leq n^{\varsigma}}\right|\leq {1\over
n}E\|\xi_k\|^2\leq {E(\|\xi_k\|^\kappa\vee 1)\over n}.\ee

At last, \be\label{452}  \left|E{\xi_{kr}\over \sqrt {n}}\cdot
\1_{\|\zeta_{k}\|\leq n^{\varsigma}}\right|=\left|E{\xi_{kr}\over
\sqrt {n}}\cdot \1_{\|\zeta_{k}\|> n^{\varsigma}}\right|\leq
n^{-{1\over 2}} \Big[E|\xi_{kr}|^\kappa\Big]^{1\over \kappa}
\Big[P(\|\zeta_{kr}\|\geq n^{\varsigma} )\Big]^{\kappa-1\over
\kappa}\leq
$$
$$
\leq\Cs(U,\mu_\kappa(\xi))\cdot n^{-{1\over 2}}\cdot
\Big[n^{-\varsigma \kappa}\Big]^{\kappa-1\over \kappa}\leq
{\Cs(U,\mu_\kappa(\xi))\over n}
 \ee
(recall that $E\xi_{kr}=0$). The triple $(\eta_k,\eps_k,\zeta_k)$
is independent of $\Ff_{k-1}$. Thus, taking in (\ref{46})
conditional expectation w.r.t. $\Ff_{k-1}$ and taking into account
inequalities (\ref{45}), we obtain an estimate
$$
 E\Upsilon_k^q\leq
E\left(\Upsilon_{k-1}+\Theta_{k-1}\cdot {1\over
n}\right)^q\!\!+\Cs(a,b,d,U,\mu_\kappa(\xi),q){E\Upsilon_{k-1}^q\over
n}\leq \left(1+{\Cs_1(a,b,d,U,\mu_\kappa(\xi),q)\over
n}\right)E\Upsilon_{k-1}^q+
$$
\be\label{47}
+\Cs_2(a,b,d,U,\mu_\kappa(\xi),q)\sum_{l=0}^{q-1}\left(1+{1\over
n}\right)^l\cdot{1\over n^{q-l}}\cdot E\Upsilon_{k-1}^l. \ee Let
us show that (\ref{47}) provide the family of estimates
\be\label{48} E\Upsilon_k^q\leq
\Cs(a,b,d,U,\mu_\kappa(\xi),q)\brakes{k\over n}^q, \quad
k=1,\dots, n, \quad q\in \NN \ee (note that (\ref{48}) is exactly
(\ref{43}) with $m=1$ and $p=2q$). We use induction by $q$. For
$q=1$, (\ref{47}) implies that
$$
E\Upsilon_k\leq {\Cs_2(a,b,d,U,\mu_\kappa(\xi),1)\over
n}+{\Cs_2(a,b,d,U,\mu_\kappa(\xi),1)\over n}\cdot
\left(1+{\Cs_1(a,b,d,U,\mu_\kappa(\xi),1)\over n}\right)+ \dots +
$$
$$
+{\Cs_2(a,b,d,U,\mu_\kappa(\xi),1)\over n}\cdot
\left(1+{\Cs_1(a,b,d,U,\mu_\kappa(\xi),1)\over n}\right)^{k-1}\leq
{k\over n}\cdot \Cs_2(a,b,d,U,\mu_\kappa(\xi),1)\cdot
e^{\Cs_1(a,b,d,U,\mu_\kappa(\xi),1)}.
$$
Similarly, if (\ref{48}) holds true for all $q\leq Q-1$, then
(\ref{47}) implies that
$$
E\Upsilon_k^Q\leq \Cs_2(a,b,d,U,\mu_\kappa(\xi),Q)\cdot
e^{\Cs_1(a,b,d,U,\mu_\kappa(\xi),Q)}\cdot \sum_{l=0}^{Q-1}
{2^l\over n^{Q-l}}\cdot \Cs(a,b,d,U,\mu_\kappa(\xi),l)\cdot
\brakes{k\over n}^{l}\leq
$$
$$
\leq \Cs(a,b,d,U,\mu_\kappa(\xi),Q)\brakes{k\over n}^Q,
$$
that proves (\ref{48}) for $q=Q$. This proves the statement of the
lemma for $m=1$. For arbitrary $m$, the proof is analogous: one
should write difference relations for the higher derivatives of
$X_n$, analogous to (\ref{41}), and then again use the moment
estimates of the same type with the given above. This step does
not differ principally from the one for SDE's driven by a Wiener
process (see, for instance \cite{iked_wat}, Chapter V \S8), and
thus  we omit its detailed exposition here. The lemma is proved.

\begin{cor}\label{c41}
$$
N_{j}(X_n(t))\leq \Cs(a,b,d,U,\mu_\kappa(\xi),j)\cdot \sqrt{t},
\quad t\in\left[{1\over n}, 1\right], j\in\NN.
$$
\end{cor}

Let us proceed with the investigation of the Malliavin's matrix
$\sigma^{f^n}$ for $f^n=X_n(t)$. Denote by $\Ef_{i,j}^n, 0\leq
i\leq j\leq n$ the difference analogue of the stochastic exponent
for (\ref{41}), i.e., the family of $\Re^{d\times d}$-valued
variables satisfying the relations \be\label{49}
\Ef_{i,i}^n=I_{\Re^d},\quad \Ef_{i,j}^n=\Ef_{i,j-1}^n+\nabla
a\brakes{X_n\brakes{j-1\over n}}\Ef_{i,j-1}^n\cdot {1\over n}
 +\sum_{r=1}^d \left[\nabla
b_r\brakes{X_n\brakes{j-1\over n}}\Ef_{i,j-1}^n\cdot
{\xi_{jr}\over\sqrt n}\right],  \ee $j=i,\dots, n$. Then one can
easily obtain  the representation for $Y_n(\cdot)$, \be\label{410}
Y_{n}\brakes{k\over
n}=\theta_n(\zeta)\sum_{j=1}^k\sum_{r=1}^d{\1_{\eps_k=1}\psi(\eta_k)\over\sqrt
n}\cdot\left[\Ef_{j,k}^nb_r\brakes{X_n\brakes{j-1\over
n}}\right]\otimes e_{kr}.\ee Denote $f^{n,k}=X\brakes{k\over n}$
and $\sigma_{n,k}=\sigma^{f^{n,k}}$. By (\ref{410}), we have
$$\sigma_{n,k}={\theta_n(\zeta)\over n}\sum_{j=1}^k\sum_{r=1}^d
[\psi^2(\eta_j)\1_{\eps_j=1}]
\cdot\left[\Ef_{j,k}^nb_r\brakes{X_n\brakes{j-1\over
n}}\right]\otimes\left[\Ef_{j,k}^nb_r\brakes{X_n\brakes{j-1\over
n}}\right]=
$$
$$ ={\theta_n(\zeta)\over n}\sum_{j=1}^k
[\psi^2(\eta_j)\1_{\eps_j=1}]
\cdot\left[\Ef_{j,k}^nb\brakes{X_n\brakes{j-1\over
n}}\right]\cdot\left[\Ef_{j,k}^nb\brakes{X_n\brakes{j-1\over
n}}\right]^*=
$$
$$
={\theta_n(\zeta)\over n}\sum_{j=1}^k
[\psi^2(\eta_k)\1_{\eps_k=1}] \cdot\Ef_{j,k}^n\cdot
[bb^*]\brakes{X_n\brakes{j-1\over n}}\cdot [\Ef_{j,k}^n]^*.$$
Together with the family $\{\Ef_{i,j}^n\}$, we consider the family
$\{\tilde\Ef_{i,j}^n\}$ defined by \be\label{411}\tilde
\Ef_{i,i}^n=I_{\Re^d},\quad
\tilde\Ef_{i,j}^n=\tilde\Ef_{i,j-1}^n+\nabla
a\brakes{X_n\brakes{j-1\over n}}\tilde\Ef_{i,j-1}^n\cdot {1\over
n}
 +\sum_{r=1}^d \left[\nabla
b_r\brakes{X_n\brakes{j-1\over n}}\tilde\Ef_{i,j-1}^n\cdot
{\xi_{jr}^n\over\sqrt n}\right], \ee $j=i,\dots, n$, where
$\xi_{i}^n=\xi_i\1_{\|\zeta_i\|\leq n^\varsigma}, i=1,\dots, n.$
By the construction, $\tilde \Ef_{i,j}^n=\Ef_{i,j}^n$ on the set
$\{\theta_n(\zeta)=1\}$, therefore $$
\sigma_{n,k}={\theta_n(\zeta)\over n}\sum_{j=1}^k
[\psi^2(\eta_j)\1_{\eps_j=1}] \cdot\tilde\Ef_{j,k}^n\cdot
[bb^*]\brakes{X_n\brakes{j-1\over n}}\cdot
[\tilde\Ef_{j,k}^n]^*.$$  We have
$$
\tilde \Ef_{i,j}^n=\prod_{l=j}^{i+1}\left[I_{\Re^d}+\nabla
a\brakes{X_n\brakes{l-1\over n}}\cdot {1\over n}+\sum_{r=1}^d
\nabla b_r\brakes{X_n\brakes{l-1\over n}}\cdot
{\xi_{lr}^n\over\sqrt n}\right].
$$
Since $\nabla a, \nabla b$ are bounded and
$|\xi_{jr}^n|\1_{\|\zeta_j\|\leq n^\varsigma}\leq \max(\max_{x\in
U}\|x\|, n^\varsigma)$, there exists $n_0=n_0(a,b,d,U,\varsigma)$
such that
$$
\left\|\nabla a\brakes{X_n\brakes{l-1\over n}}\cdot {1\over
n}+\sum_{r=1}^d \nabla b_r\brakes{X_n\brakes{l-1\over n}}\cdot
{\xi_{lr}^n\over\sqrt n}\right\|\1_{\|\zeta_j\|\leq
n^\varsigma}\leq {1\over 2}, \quad n\geq n_0.
$$
Then $\tilde\Ef_{i,j}^n$ is invertible and
$$
[\tilde \Ef_{i,j}^n]^{-1}=\prod_{l=i+1}^{j}\left[I_{\Re^d}+\nabla
a\brakes{X_n\brakes{l-1\over n}}\cdot {1\over n}+\sum_{r=1}^d
\nabla b_r\brakes{X_n\brakes{l-1\over n}}\cdot
{\xi_{lr}^n\over\sqrt n}\right]^{-1}.
$$
Thus, on the set
$\{\theta_n(\zeta)=1\}\backslash\{\eps_1=\dots=\eps_n=0\}$, the
matrix $\sigma_{n,k}$ is invertible and \be\label{4111}
\|\sigma_{n,k}^{-1}\|=[\inf_{\|v\|=1} (\sigma_{n,k}v,v)]^{-1}\leq
[\max_{i\leq j\leq n}\|[\tilde \Ef_{i,j}^n]^{-1}\|]^2\cdot
\left\{{1\over n}\sum_{j=1}^k
\psi^2(\eta_j)\1_{\eps_j=1}\right\}^{-1}.\ee

\begin{lem}\label{l43} For every $p\geq 1$,
$$
E[\max_{i\leq j\leq n}\|[\tilde \Ef_{i,j}^n]^{-1}\|]^p\leq
\Cs(a,b,d,U,\mu_\kappa(\xi),p), \quad n\in\NN.
$$
\end{lem}
\demo Since $[\tilde \Ef_{i,j}^n]^{-1}=\tilde
\Ef_{0,i}^n\cdot[\tilde \Ef_{0,j}^n]^{-1}$, it is enough to prove
that \be\label{412}E[\max_{i\leq n}\|\tilde \Ef_{0,i}^n\|]^p\leq
\Cs(a,b,d,U,\mu_\kappa(\xi),p),\ee \be\label{413}E\max_{j\leq
n}|\det\tilde \Ef_{0,j}^n|^{-p}\leq
\Cs(a,b,d,U,\mu_\kappa(\xi),p).\ee Let us prove inequality
 \be\label{4121}
E[\max_{i\leq n}\|\tilde \Ef_{0,i}^n\|_2]^p\leq
\Cs(a,b,d,U,\mu_\kappa(\xi),p) \ee
 with
$\|A\|_2\equiv \sqrt{\sum_{lr}A_{lr}^2}$; this will provide
inequality (\ref{412}). We deduce from (\ref{411}) that
$Z^n_i\equiv\|\tilde \Ef^n_{0,i}\|^2_2$ satisfy relations
analogous to (\ref{44}), i.e. \be\label{414}
Z_i^n=Z_{i-1}^n+V_{i-1}^{1,n}\cdot {1\over
n}+\sum_rV_{i-1}^{2,n}\cdot {\xi_{ir}^n\over
\sqrt{n}}+\sum_{r_1,r_2}V_{i-1}^{3,n}\cdot
{\xi_{ir_1}^n\xi_{ir_1}^n\over {n}},\quad  i=1,\dots,n \ee with an
$\{\Ff_i\}$ -- adapted sequences $V_i^{1,n}, V_{i,\cdot}^{2,n},
V_{i,\cdot,\cdot}^{3,n}$ such that \be\label{4141} |V_i^{1,n}|\leq
\Cs(a,b,d,U)(1+Z_i^n),\quad |V_{i,r}^{2,n}|\leq \Cs(a,b,d,U)
Z_i^n,\quad |V_{i,r_1,r_2}^{3,n}|\leq \Cs(a,b,d,U) Z_i^n. \ee Then
then the moment estimates analogous to those made in the proof of
Lemma \ref{l42} provide that \be\label{415} \max_{i\leq
n}(EZ_i^n)^{p\over 2}\leq \Cs(a,b,d,U,\mu_\kappa(\xi),p).\ee
Denote $A_j^n=Z_0^n+\sum_{i=1}^j\Delta A_i^n,
M_j^n=\sum_{i=1}^j\Delta M_i^n$ with
$$
\Delta M_i^n=\sum_rV_{i-1}^{2,n}\cdot {\xi_{ir}^n-E\xi_{ir}^n\over
\sqrt{n}}+\sum_{r_1,r_2}V_{i-1}^{3,n}\cdot
{\xi_{ir_1}^n\xi_{ir_2}^n-E\xi_{ir_1}^n\xi_{ir_2}^n\over {n}}
$$
and
$$
\Delta A_i^n=Z_i^n-Z_{i-1}^n-M_i^n=V_{i-1}^{1,n}\cdot {1\over
n}+\sum_rV_{i-1}^{2,n}\cdot {E\xi_{ir}^n\over
\sqrt{n}}+\sum_{r_1,r_2}V_{i-1}^{3,n}\cdot
{E\xi_{ir_1}^n\xi_{ir_2}^n\over {n}}.
$$
Then $Z_i=A_i+M_i$. By (\ref{452}) and (\ref{4141}),
$$
|\Delta A_i^n|\leq {\Cs(a,b,U,\mu_\kappa(\xi))\over n}\cdot
(1+Z_i),
$$
and therefore (\ref{415}) provides that
$$
E\max_{j\leq n} |A_j^n|^{r\over 2}\leq
\Cs(a,b,U,\mu_\kappa(\xi),p).
$$
Similarly, Burkholder inequality together with (\ref{4141}) and
(\ref{450}) -- (\ref{452}) provides that
$$
E\max_{j\leq n} |M_j^n|^{p\over 2}\leq
\Cs(a,b,U,\mu_\kappa(\xi),p),
$$
that proves (\ref{4121}), and therefore (\ref{412}).

On the  set $\{\|A\|_2\leq {1\over 2}\}\subset \Re^{d\times d}$,
the function $\Phi:A\mapsto (\det[I_{\Re^d}+A])^{-1}$ can be
represented in the form
$$
\Phi(A)=1+Q(A)+\vartheta(A),
$$
where $P$ is a polynomial of $A$ with $\deg A\leq \kappa d$ and
$|\vartheta(A)|\leq \Cs \|A\|^{\kappa+1}_2$. We have
$$\left\|\nabla a\brakes{X_n\brakes{l-1\over n}}\cdot {1\over
n}+\sum_{r=1}^d \nabla b_r\brakes{X_n\brakes{l-1\over n}}\cdot
{\xi_{lr}^n\over\sqrt n}\right\|^{\kappa+1}_2\leq
{\Cs(a,b,d,U)\over n}.
$$
Therefore,
$$(\det\tilde \Ef_{0,j}^n)^{-1}=(\det \tilde \Ef_{0,{j-1}}^n)^{-1}\cdot\left[1+
 Q_{j-1}^n\left({\xi_{j1}^n\over \sqrt n},\dots,{\xi_{jd}^n\over
\sqrt n}\right)+\vartheta_i^n\right],$$ where $|\vartheta_i^n|\leq
{\Cs(a,b,d,U)\over n}$,  $\vartheta_i^n$ is $\Ff_j$ -- measurable,
$Q_{j-1}^n$ is a polynomial with $\deg Q_{j-1}^n\leq \kappa d$ and
its coefficients are $\Ff_{j-1}$ -- measurable and bounded by some
constant depending on the coefficients $a,b$. Repeating the
arguments used in the proof of (\ref{412}) we obtain (\ref{413}).
The lemma is proved.

\begin{lem}\label{l44} For every $p\in\NN, c>0$,
$$E\left\{\sum_{j=1}^k
\psi^2(\eta_j)\1_{\eps_j=1}\right\}^{-p}\1_{\sum_{j=1}^k\eps_j\geq
ck}\leq \Cs(c,\psi,p)k^{-p},\quad k\geq {2p+1\over c}.
$$
\end{lem}

\emph{Remark.}  For arbitrary $\psi\in C^\infty(\Re^d)$ with
$\psi=0$ on $\partial U$, the given above statement  may fail. It
is crucial for $\psi$ to have non-zero normal derivative  at (some
part of) the boundary in order to provide (\ref{4151}).

\demo Since $\eta$ and $\eps$ are independent,
$$E\left\{\sum_{j=1}^k
\psi^2(\eta_j)\1_{\eps_j=1}\right\}^{-p}\1_{\sum_{j=1}^k\eps_j\geq
ck}\leq E\left\{\sum_{j=1}^{]ck[} \psi^2(\eta_j)\right\}^{-p},
$$
where $]x[\eqdef \min\{n\in \ZZ|n\geq x\}$. By the construction of
the function $\psi$, \be\label{4151} P(\psi^2(\eta_j)\leq z)\sim
\Cs(\psi)\cdot \sqrt{z},\quad z\to 0+. \ee Therefore
$$
P(\psi^2(\eta_1)+\dots+\psi^2(\eta_l)\leq z)\sim \Cs(\psi,l)\cdot
z^{l\over 2},\quad z\to 0+
$$
and \be\label{416}
E(\psi^2(\eta_1)+\dots+\psi^2(\eta_l))^{-p}<+\infty \ee as soon as
$l>2p$. We put $q=2p+1, N=\left[{]ck[\over q}\right]$ and divide
the set $\{1,\dots,]ck[\}$ on the blocks
$$
\{1,\dots,
q\},\,\{q+1,\dots,2q\},\dots,\{(N-1)q+1,\dots,Nq\},\,\{Nq+1,\dots,]ck[\}
$$
(the last block may be empty). We denote
$$
\vartheta_i=\sum_{j=(i-1)q+1}^{iq}\psi^2(\eta_j),\quad
i=1,\dots,N.
$$
We have
$$
E\left\{\sum_{j=1}^{]ck[} \psi^2(\eta_j)\right\}^{-p}\leq
E\left(\sum_{i=1}^N \vartheta_i\right)^{-p}=N^{-p}\cdot
E\left({1\over N}\sum_{i=1}^N \vartheta_i\right)^{-p}.
$$
The function $x\mapsto x^{-p}$ is convex on $\Re^+$, and therefore
$$
E\left({1\over N}\sum_{i=1}^N \vartheta_i\right)^{-p}\leq
E\left({1\over N}\sum_{i=1}^N
\vartheta_i^{-p}\right)=E\vartheta_1^{-p}<+\infty
$$
(the last inequality follows from (\ref{416})). If $k\geq {2p+1\over
c}$, then ${]ck[\over q}\geq 1$ and therefore $N=\left[{]ck[\over
q}\right]\geq {]ck[\over 2q}\geq {c\over 4p+2}\cdot k$. Thus
$$
E\left\{\sum_{j=1}^{]ck[} \psi^2(\eta_j)\right\}^{-p}\leq
\Cs(\psi,p)\cdot \left({c\over 4p+2}\right)^{-p}\cdot k^{-p}.
$$
The lemma is proved.

Inequality (\ref{4111}), Lemma \ref{l43} and Lemma \ref{l44}
provide the following estimate. For a given $c>0$ and $t\in[0,1]$,
we put $\tilde\Xi_n=\{\theta_n(\zeta)=1\}\bigcap
\left\{\sum_{j=1}^{[tn]}\eps_j\geq c[tn]\right\}$.
\begin{cor}\label{c42} For $p\in \NN$ and $[tn]> {2p+1\over c}$,
$$
E\|\varrho^{f^{n},\tilde\Xi_n}\|_{\mathcal{M}}^{p}\leq
\Cs(a,b,d,U,\psi,\mu_\kappa(\xi),p) \cdot t^{-p}.
$$
\end{cor}

At last, let us give an estimates for the tail probabilities for
$f^n$. The following lemma is completely analogous to Lemma
\ref{l42}; the proof is omitted.

\begin{lem}\label{l421} For every $p\geq 1$, there exist constant
 $\Cs(a,b,d,U,\mu_\kappa(\xi),p)$ such that
$$E\|X_n(t)-X(0)\|^p_{\Re^d}\1_{\{\theta_n(\zeta)=1\}}\leq
 \Cs(a,b,d,U,\mu_\kappa(\xi),p)\cdot t^{{p\over 2}}, \quad t\in \left[{1\over n}, 1\right].
$$
  \end{lem}
\begin{cor}\label{c411} For every $p\geq 1$, there exists constant
$C_p$, dependent on $a,b,d,U,\mu_\kappa(\xi),p$, such that
$$
P(\|X_n(t)-X_n(0)\|\geq y,\theta_n(\zeta)=1)\leq
C_p\left(1+{\|y\|^2\over t}\right)^{-p}, \quad n\in \NN, t\in
\left[{1\over n},1\right].
$$
\end{cor}

\begin{lem}\label{l45} There exist constants $\Cs_1,\Cs_2, \Cs_3$,
dependent on $a,b,U,d,\mu_\kappa(\xi)$,  such that, for every
$\lambda\in\Re^d$ with $\|\lambda\|\leq \Cs_1 n^{1\over
\kappa+1}$, \be\label{4161} Ee^{(\lambda,
X_n(t)-X_n(0))}\1_{\{\theta_n(\zeta)=1\}}\leq \Cs_2 e^{\Cs_3 t
\|\lambda\|^2 },\quad t\in\left[{1\over n},1\right], n\in \NN. \ee
\end{lem}
\demo For a given $\lambda$, denote $Z_n(t)=e^{(\lambda,
X_n(t)-X_n(0))}$. We have $Z_n(0)=1$. On the other hand,
$$
\left|{\xi_{kr}\over \sqrt n}\right|\leq {\max_{y\in U}\|y\|\over
\sqrt n} + n^\varsigma\cdot n^{-{1\over 2}}={\max_{y\in
U}\|y\|\over \sqrt n}+n^{-{1\over \kappa+1}}
$$
on the set $\{\theta_n(\zeta)=1\}$. Thus there exists a constant
$\Cs_4$ such that, for $\|\lambda\|\leq \Cs_1n^{1\over \kappa+1}$,
$$
\left|{1\over n}\left(\lambda, a\brakes{X_n\left(k-1\over
n\right)}\right)+\sum_{r=1}^d {\xi_{kr}\over \sqrt
n}\left(\lambda,b_r\brakes{X_n\left(k-1\over
n\right)}\right)\right|\leq \Cs_4
$$
on the set $\{\theta_n(\zeta)=1\}$. Using the elementary
inequality $e^x\leq 1+x+\Cs x^2, |x|\leq \Cs_4$, we obtain that,
on the same set,
$$
Z_n\left(k\over n\right)=Z_n\left(k-1\over
n\right)\exp\left[{1\over n}\left(\lambda,
a\brakes{X_n\left(k-1\over n\right)}\right)+\sum_{r=1}^d
{\xi_{kr}\over \sqrt n}\left(\lambda,b_r\brakes{X_n\left(k-1\over
n\right)}\right)\right]\leq
$$
$$
\leq Z_n\left(k-1\over n\right)\left[1+\Theta_{k-1}\cdot
{\|\lambda\|+\|\lambda\|^2\over n} +
\sum_{r=1}^d\Lambda_{k-1,r}\cdot
{\xi_{kr}(\|\lambda\|+\|\lambda\|^2)\over\sqrt n}+\right.
$$
\be\label{417} +\left.\sum_{r_1,r_2=1}^d\Delta_{k-1,r_1,r_2}\cdot
{\xi_{kr_1}\xi_{kr_2}\|\lambda\|^2\over n}\right] \ee with an
$\Ff_{k-1}$ -- measurable coefficients $\Theta_{k-1},
\Lambda_{k-1,r_1,r_2}$, bounded by some constant $\Cs$. An
arguments, analogous to those used in the proof of Lemma
\ref{l42}, provide that (\ref{417}) implies the estimate
\be\label{418} EZ_n\left(k\over n\right)\leq
\left(1+\Cs{\|\lambda\|+\|\lambda\|^2\over n}\right)^k \leq
\exp\left[2\Cs\cdot {k\over n}\cdot (1+\lambda^2)\right], k,n \in
\NN .\ee This is exactly (\ref{4161}) for $t={k\over n}$. For
$t\in \left({k\over n}, {k+1\over n}\right)$, $Z_n(t)$ is a linear
combination of $Z_n\left(k\over n\right)$ and $Z_n\left(k+1\over
n\right)$. Therefore,  (\ref{4161}) follows from (\ref{418}) via
Jensen's inequality and relation ${k+1\over n}\leq 2{k\over n}\leq
2t$ (recall that $t\geq {1\over n}$ and thus $k\geq 1$). The lemma
is proved.

\begin{cor}\label{c43} There exist constants $\Cs_5,\Cs_6, \Cs_7$,
dependent on $a,b,U,d,\mu_\kappa(\xi)$,  such that $$
P(\|X_n(t)-X_n(0)\|\geq y,\theta_n(\zeta)=1)\leq \Cs_6e^{-\Cs_7
{y^2\over t}},\quad y\in (0, \Cs_5 t n^{1\over \kappa+1}),n\in
\NN, t\in \left[{1\over n},1\right]
$$
and
$$
P(\|X_n(t)-X_n(0)\|\geq y,\theta_n(\zeta)=1)\leq \Cs_6e^{-\Cs_7
n^{1\over \kappa+1} y},\quad y\geq  \Cs_5 t n^{1\over
\kappa+1},n\in \NN, t\in \left[{1\over n},1\right].
$$
\end{cor}

\demo It is enough to verify that, for any coordinate $(X_n)_j$ of
the process $X_n$, $j=1,\dots, d$, there exist constants $\tilde
\Cs_5,\tilde \Cs_6, \tilde \Cs_7,\tilde \Cs_8$ such that
\be\label{419} P(\pm((X_n)_j(t)-(X_n)_j(0))\geq
y,\theta_n(\zeta)=1)\leq \tilde \Cs_6e^{-\tilde \Cs_7 {y^2\over
t}},\quad y\in (0, \tilde\Cs_5 t n^{1\over \kappa+1}), n\in \NN,
t\in \left[{1\over n},1\right] \ee and \be\label{420}
P(\pm((X_n)_j(t)-(X_n)_j(0))\geq y,\theta_n(\zeta)=1)\leq \tilde
\Cs_6e^{-\tilde\Cs_8 n^{1\over \kappa+1} y},\quad y\geq
\tilde\Cs_5 t n^{1\over \kappa+1},n\in \NN, t\in \left[{1\over
n},1\right]. \ee Inequality (\ref{419}) with $\tilde
\Cs_5=2\Cs_1\Cs_3$ and $\tilde \Cs_7=[2 \Cs_3]^{-1}$ follows from
(\ref{4161}) with $\lambda=\left(\pm{y\over 2\Cs_3 t}\right)\cdot
e_j, $ where $e_j$ is the $j$-th coordinate vector in $\Re^d$.
Inequality (\ref{420}) with the same $\tilde \Cs_5$ and $\tilde
\Cs_8={\Cs_1\over 2}$ follows from (\ref{4161}) with
$\lambda=\left(\pm\Cs_1 n^{1\over \kappa+1}\right)\cdot e_j.$

\emph{Proof of Theorem \ref{t11}.} We take $p=8(d+1)$ and fix some
$c\in (0,\alpha)$ ($\alpha$ is given in condition (B3)). We write
$n_*=n_0(a,b,U,\varsigma)$ (see the notation before Lemma
\ref{l43}) and put
$$
\Xi_n^t=\begin{cases}\{\theta_n(\zeta)=1\}\bigcap
\left\{\sum_{j=1}^{[tn]}\eps_j\geq c[tn]\right\},& n\geq n_*,[tn]>
{2p+1\over c}\\
\emptyset,&\hbox{otherwise}\end{cases},
$$
$$ Q_{x,t}^n(dy)=P(X_n(x,t)\in dy, \Xi_n^t)=P_{\Xi_n^t}(f_{x,t}^n\in dy-x),
$$
$$
R_{x,t}^n(dy)=P(X_n(x,t)\in dy, \Omega\backslash\Xi_n^t).
$$

Corollaries \ref{c41} and \ref{c42} provide  condition (C1) of the
Theorem \ref{t31}. By the statement (a) of this theorem,
\be\label{421} Q_{x,t}^n(dy)=q_{x,t}^n(y)\, dy\hbox{ with
}q_{x,t}^n\leq \Ks_{d}(f^n,\Xi_n)\cdot P^{1\over 2}(\|f^n\|\geq
\|y\|, \theta_n(\zeta)=1). \ee Moreover, Corollaries \ref{c41} and
\ref{c42} provide an explicit estimate for $\Ks_{d}(f^n,\Xi_n)$.
Namely, for some constant $\Cs$ dependent on
$a,b,c,d,\mu_\kappa(\xi),U,\psi$, \be\label{4211}
\Ks_{d}(f^n,\Xi_n)\leq \Cs\sum_{M=0}^d[\sqrt t]^{d+2M}\cdot
[t^{-4(M+d)}]^{1\over 4}=(d+1)\Cs\, t^{-{d\over 2}}. \ee Thus the
statement (a) of Theorem \ref{t31} and Corollaries
\ref{c43},\ref{c411} provide statement (ii) of Theorem \ref{t11}.

By  Chebyshev inequality,
$$
P(\theta_n(\zeta)=0)\leq n\cdot
n^{-\kappa\cdot\varsigma}=n^{-\epsilon(\kappa)}.
$$
Take $\lambda= \ln\left({\alpha(1-c)\over c(1-\alpha)}\right)>0$.
By  Chebyshev inequality, we get, after some simple calculations,
\be\label{4212} P\left(\sum_{j=1}^{k}\eps_j< ck\right)\leq {E
e^{-\lambda \sum_{j=1}^k\eps_j}\over e^{-\lambda c
k}}=[\Psi(\alpha,c)]^k, \quad k\in \NN, \ee with
$\Psi(\alpha,c)=\left({1-\alpha\over
1-c}\right)^{1-c}\left({\alpha\over c}\right)^{c}$.  One can
verify that $\Psi(\alpha,c)<0$ for $0<c<\alpha<1$. Thus, we can
conclude that, for $\rho=-{1\over 2}\ln \Psi(\alpha,c)>0,$
$$P(\Omega\backslash \Xi_n^t)\leq n^{-\epsilon(\kappa)}+e^{-\rho nt} \quad\hbox{when}\quad n\geq n_*,[tn]> {2p+1\over
c}$$ (we have used here that $[nt]\geq {nt\over 2}$ for $t\geq
{1\over n}$). Thus, for all $t>0$,
\be\label{422}P(\Omega\backslash \Xi_n^t)\leq
D[n^{-\epsilon(\kappa)}+e^{-\rho nt}] \ee with the constant $D$
dependent on $n_*, p,c,\rho$. This provides statement (iii) of
Theorem \ref{t11}.

We have shown that if $x\in\Re^d, t>0$ are fixed then the
functions $f^n=f_{x,t}^n$ and the sets $\Xi_n=\Xi_n^t$ satisfy all
the conditions of Theorem \ref{t31}. This means that
$q_{x,t}^n(y)\to p_{x,t}(y)$ uniformly w.r.t. $y\in \Re^d$. In
order to show that this convergence holds uniformly w.r.t. $t\in
[\delta,1], x,y\in\Re^d$, we need to show that, for every
sequences $\{t_n\}\subset [\delta,1],\{x_n\},\{y_n\}\subset
\Re^d$, \be\label{423} q_{x_n,t_n}^n(y_n)-p_{x_n,t_n}(y_n)\to 0,
\quad n\to \infty. \ee We can suppose that $\{t_n\}$ converges to
some $t\in [\delta,1]$. The functions $a,b$ are bounded together
with their derivatives up to the second order, and therefore the
sequences of the functions
$$
a_n(\cdot)=a(\cdot+x_n),\quad b_n(\cdot)=b(\cdot+x_n)
$$
are pre-compact in $C^1(\Re^d,\Re^d)$ and $C^1(\Re^d,\Re^{d\times
d})$, correspondingly. We can suppose that
$$
a_n\to \tilde a\hbox{ in }C^1(\Re^d,\Re^d),\quad b_n\to \tilde
b\hbox{ in }C^1(\Re^d,\Re^{d\times d}).
$$
Consider the processes $Z_n$ defined by the relations of the type
(\ref{01}),(\ref{02}) with $Z_n(0)=0$ and the coefficients $a,b$
replaced by $a_n, b_n$. Also, consider the processes $Z^n$ defined
by the SDE's of the type (\ref{00}) with $Z^n(0)=0$ and the
coefficients $a,b$ replaced by $a_n, b_n$. At last, consider the
process $Z$ defined by the SDE's of the type (\ref{00}) with
$Z(0)=0$ and the coefficients $a,b$ replaced by $\tilde a, \tilde
b$. Denote $f_n=Z_n(t_n), \Xi_n=\Xi_n^{t_n}.$ It is easy to verify
that $Z_n$ converge weakly in $C([0,1],\Re^d)$ to $Z$ (see, for
instance, Proposition 5.1 \cite{kurtz_protter}). Thus, for the
sequences $f_n,\Xi_n$, all the conditions of Theorem \ref{t31}
hold true with $f=Z(t)$. This means that $f$ possesses a
distribution density $p^f$ and \be\label{424} \sup_y
|p_{\Xi_n}^{f_n}(y)-p^f(y)|\to 0. \ee Similarly,  one can show
that, for $f^n=Z^n(t_n)$, the distribution density $p^{f^n}$
exists and \be\label{425} \sup_y |p^{f^n}(y)-p^f(y)|\to 0. \ee
Now, (\ref{423}) is provided by the relations
(\ref{424}),(\ref{425}) and
$$
q_{x_n,t_n}(y_n)=p_{\Xi_n}^{f_n}(y_n-x_n),\quad
p_{x_n,t_n}(y_n)=p^{f^n}(y_n-x_n).
$$
This proves statement (i) of Theorem \ref{t11}.

The proof of (ii$'$) and (iii$'$) can be conducted analogously,
with an appropriate changes of the truncation procedure and
corresponding estimates. Under (B2$_{\exp}$), we put, instead of
(\ref{40}),
$$
\theta_n(\zeta)=\1_{\max_{k\leq n}\|\zeta_k\|\leq \delta \sqrt n}.
$$
By the Chebyshev's inequality,
$$
P(\theta_n(\zeta)=0)\leq n\cdot e^{-\kap\cdot (\delta\sqrt
n)^2}\leq e^{-\tilde \rho n}.
$$
with an appropriate $\tilde \rho>0$. This and the estimate
(\ref{4212}) provide statement (iii$'$).

Under (B2$_{\exp}$) and the truncation level given above, the
estimates (\ref{451}) -- (\ref{452}) have their (simpler)
analogues, and thus the statement of Lemma \ref{l42} holds. The
constant $\delta$ in the definition of the truncation level can be
made small enough to provide inequality
$$
\left\|\nabla a\brakes{X_n\brakes{l-1\over n}}\cdot {1\over
n}+\sum_{r=1}^d \nabla b_r\brakes{X_n\brakes{l-1\over n}}\cdot
{\xi_{lr}^n\over\sqrt n}\right\|\1_{\|\zeta_j\|\leq \delta \sqrt
n}\leq {1\over 2}, \quad n\geq n_0
$$
to hold with some $n_0$ dependent on $a,b,U$. Then the statement
of Lemma \ref{l43} holds true, also. Lemma \ref{l44} does not
depend on the truncation procedure. Thus the Corollaries
\ref{c41},\ref{c42} hold true and provide the principal estimate
(\ref{4211}).

Under condition (B2$_{\exp}$),
$$
E\left[\exp\left[{1\over n}\left(\lambda,
a\brakes{X_n\left(k-1\over n\right)}\right)+\sum_{r=1}^d
{\xi_{kr}\over \sqrt n}\left(\lambda,b_r\brakes{X_n\left(k-1\over
n\right)}\right)\right]\Big|\Ff_{k-1\over n}\right]\leq \tilde
\Cs_1 e^{\tilde Cs_2{\|\lambda\|^2\over n}}\quad \mathrm{a.s.}
$$
for every $\lambda\in\Re^d$ with some constants $\tilde C_1,\tilde
C_2$ dependent on $a,b,d,\kap,Ee^{\kap \|\xi_k\|^2}$. Then the
arguments analogous to those made in the proof of Lemma \ref{l42}
provide that the estimate (\ref{4161}) holds true for every
$\lambda\in \Re^d$. Consequently, the first inequality in
Corollary \ref{c43} holds true for every $y>0$. This inequality,
the estimate (\ref{4211}) and Theorem \ref{t31} provide (iii$'$).
This completes the proof of Theorem \ref{t11}.

\subsection{Proof of Theorem \ref{t13}} The implication
2 $\Rightarrow$ 1 is obvious. Let us first  prove 2 under
additional supposition (B3). We put, in the notation of Section 3,
$\theta_n(\zeta)\equiv 1, f^n={1\over \sqrt n}\sum_{k=1}^n\xi_k,$
$$
\Xi_n=\begin{cases} \left\{\sum_{j=1}^{n}\eps_j\geq cn\right\},&
n\geq n_*\\
\emptyset,&\hbox{otherwise}\end{cases}
$$
with $n_*,c$ that will be defined later. Then
$$
Df^n={1\over \sqrt
n}\sum_{k=1}^n\sum_{r=1}^d\psi(\eta_k)\1_{\eps_k=1}b_r\otimes
e_{kr},
$$
where $b_r$ stands for the $r$-th coordinate vector in $\Re^d$
(the proof is straightforward and omitted). Since $\psi$ is
bounded on $U$ together with all its derivatives, this provides
the estimates analogous to those given in  Lemma \ref{l42}. The
statement of Lemma \ref{l43} is trivial now, since
$\Ef_{i,j}=I_{\Re^d}$ (the identity matrix in $\Re^d$) for every
$i,j$. Now, take $p=8(d+1), c\in(0,\alpha), n_*>{2p+1\over c}$.
Using Lemma \ref{l44}, we obtain the estimate (\ref{4211}) with
the constant $\Cs$ dependent on $c,d,\psi$. The estimate
(\ref{4212}) provides that $P(\Xi_n)\leq De^{-\rho n}$ for an
appropriate $D,\rho>0$. Now the statement 2 follows from Theorem
\ref{t31}.

Let us replace the additional supposition (B3) by the condition 1.
It is enough to prove the statement 2 for $n\in m\NN$ with some
given $m\in \NN$. We take $m=2n_0$ with $n_0$  given in the
statement 1 and have
$$
{d[P_{m}]^{ac}\over d\lambda^d}\geq {d[P_{n_0}]^{ac}\over
d\lambda^d}*{d[P_{n_0}]^{ac}\over d\lambda^d}.
$$
The function ${d[P_{n_0}]^{ac}\over
d\lambda^d}*{d[P_{n_0}]^{ac}\over d\lambda^d}$ is continuous since
$\lambda^d$-almost all points of $\Re^d$ are a Lebesgue points
(i.e., a points of $\lambda^d$-almost continuity) for any function
$g\in L_1(\Re^d)$. In addition, ${d[P_{n_0}]^{ac}\over
d\lambda^d}*{d[P_{n_0}]^{ac}\over d\lambda^d}$ is not an identical
zero due to the statement 1. Thus there exist $\tilde \alpha>0$
and an open set  $U\subset \Re^d$ such that ${d[P_{n_0}]^{ac}\over
d\lambda^d}*{d[P_{n_0}]^{ac}\over d\lambda^d}\geq \tilde
\alpha\1_U.$ Therefore the distribution of $\xi_1+\dots+\xi_{m}$
satisfies (B3). Using what we have proved before, we deduce that
the statement 2 holds for  $n\in m\NN$, and therefore for $n\in
\NN$. This completes the proof of Theorem \ref{t13}.

\subsection{Sketch of the proof of Theorem \ref{t22}} We will show
that, under conditions of Theorem \ref{t22}, the following
\emph{uniform local Doeblin condition} holds true. For two
measures $\mu_1,\mu_2$,  denote
$$
[\mu_1\wedge \mu_2](dy)=\min\left[{d\mu_1\over
d(\mu_1+\mu_2)}(y),{d\mu_2\over
d(\mu_1+\mu_2)}(y)\right](\mu_1+\mu_2)(dy).
$$

\begin{prop}\label{p41} For every ball $B$ there exists $n_B\in \NN,
T_B,\gamma_B>0$ such that \be\label{62} \Big[P_{x,T_B}\wedge
P_{x',T_B}\Big](\Re^d)\geq \gamma_B, \quad x,x'\in B \ee and, for
every $n\geq n_B$, there exists $T_B^n\in {1\over n}\ZZ_+,
T_B^n\leq T_B$ such that \be\label{63} \Big[P_{x,T_B^n}^n\wedge
P_{x',T_B^n}^n\Big](\Re^d)\geq \gamma_B, \quad x,x'\in B, \quad
n\geq n_B. \ee
\end{prop}

Once Proposition \ref{p41} is proved, one can finish the proof of
Theorem \ref{t22} following the proof of Theorem 1 in
\cite{Klo_Ver} literally. We omit this part of the discussion and
prove Proposition \ref{p41}, only.

{\it Proof of Proposition \ref{p41}.} Since (B5) implies
(B2$_\kappa$) for any $\kappa$, we can apply Theorem \ref{t11}.
One can easily see that
$$
\Big[P_{x,T_B^n}^n\wedge
P_{x',T_B^n}^n\Big](\Re^d)\geq\int_{\Re^d}\min[q_{x, t}^n(y),
q_{x',t}^n(y)]\, dy.
$$
Thus, for any sequence $t_n\to t>0$ we have, by the statement (i),
$$
\lim\inf_{n\to +\infty}\inf_{x,x'\in B}\Big[P_{x,T_B^n}^n\wedge
P_{x',T_B^n}^n\Big](\Re^d)\geq  \inf_{x,x'\in
B}\int_{\Re^d}\min[p_{x, t}(y), p_{x',t}(y)]\, dy.
$$
On the other hand (see \cite{IKO}), under condition (B1) the
function
$$
\Re^d\times(0,+\infty)\times \Re^d\ni(x,t,y)\mapsto p_{x,t}(y)
$$
is continuous and strictly positive at every point. Therefore, for
a given $B$,
$$
\gamma_B\eqdef{1\over 2}\inf_{x,x'\in B}\int_{\Re^d}\min[p_{x,
t}(y), p_{x',t}(y)]\, dy>0,
$$
and (\ref{62}), (\ref{63}) hold true for $T_B=T_B^n=1$ and
sufficiently large $n_B$. The proposition is proved.

\subsection{Sketch of the proof of Theorem \ref{t21}} It was
already mentioned in subsection 2.1 that Theorem \ref{t21} is
analogous to Theorem 2.1 \cite{kulik_da}. We refer the reader to
the paper \cite{kulik_da} for the detailed proof. Here, we expose
a principal estimate only, demonstrating that, in this proof, the
truncated local limit theorem can be used efficiently instead of
the usual one, that was used in \cite{kulik_da}.

Theorem 2.1 \cite{kulik_da} is derived from the general theorem on
convergence in distribution of a sequence of additive functionals
of Markov chains, given in the paper \cite{kar_kul} (Theorem 1).
\emph{The characteristics} of the functionals $
\phi(X),\phi_n(X_n)$ are defined by the relations
$$
f^t(x)\eqdef E[\phi^{0,t}(X_n)|X(0)=x],\quad
 f_n^{s,t}(x) \eqdef
E[\phi_{n}^{s,t}(X_n)|X_n(s)=x],\quad s={i\over n}, i\in \ZZ_+,
t>s,\quad x\in \Re^d.
$$
The first relation is due to \cite{dynkin}, Chapter 6. The second
relation was introduced in \cite{kar_kul} by an analogy with the
first one.

The key condition of Theorem 1 \cite{kar_kul} is  \be\label{61}
\underset{x\in\Re^d, s={i\over n}, t\in
(s,T)}{\sup}|f^{s,t}_n(x)-f^{t-s}(x)|\rightarrow0, \quad
n\rightarrow\infty. \ee Here, we need to verify this condition
only, since, for all the other conditions,  the proof from
\cite{kulik_da} can be used literally. We have
$$
f_{n}^{s,t}(x)={1\over n}F_n(x)+{1\over n}\sum_{k\in \NN, {k\over
n}<t-s}\int_{\Re^m}F_{n}(y) P_{x,{k\over n}}^{n}(dy)
=f_n^{0,t-s}(x), \quad s\leq t, x\in\Re^m. $$  We use the
decomposition $P^n=Q^n+R^n$ from Theorem \ref{t11} and write
$$
f_{n}^{0,t}(x)={1\over n}F_n(x)+{1\over n}\sum_{k\in \NN, {k\over
n}<t}\int_{\Re^m}F_{n}(y) R_{x,{k\over n}}^{n}(dy)+{1\over
n}\sum_{k\in \NN, {k\over n}<t}\int_{\Re^m}F_{n}(y) q_{x,{k\over
n}}^{n}(y)\, dy.
$$
The statement (ii) of Theorem \ref{t11} and the estimates,
analogous to the estimates (4.2) -- (4.10) from \cite{kulik_da},
imply that
$$
\underset{x\in\Re^d,  t\leq T}{\sup}\left|{1\over n}\sum_{k\in
\NN, {k\over n}<t}\int_{\Re^m}F_{n}(y) q_{x,{k\over n}}^{n}(y)\,
dy-f^{t}(x)\right|\rightarrow0, \quad n\rightarrow\infty.
$$
On the other hand, $\epsilon(\kappa)={8\over 7}>1$ for $\kappa=6$.
Thus, by condition (B6) and the statement (iii) of Theorem
\ref{t11},
$$
E\left[{1\over n}F_n(x)+{1\over n}\sum_{k\in \NN, {k\over
n}<t}\int_{\Re^m}F_{n}(y) R_{x,{k\over n}}^{n}(dy)\right]\leq
n^{-1}\sup_{x'}F_n(x')\left[1+\sum_{k<tn}R^n_{x,{k\over
n}}(\Re^d)\right]\leq
$$
$$
\leq n^{-1}\sup_{x'}F_n(x')\left[1+D n^{-{8\over 7}}\cdot
nt+D\sum_{k\in \NN}e^{-\gamma k}\right]\to 0, \quad n\to \infty
$$
uniformly for $x\in\Re^d, t\leq T$ for any $T\in \ax$. This proves
(\ref{61}).

\end{document}